\documentclass{article}
\usepackage{amssymb,amsthm,url}
\usepackage[colorlinks]{hyperref}
\usepackage[nameinlink]{cleveref}
\usepackage{eurosym}
\usepackage{amsmath}
\usepackage{amsfonts}
\usepackage{amssymb}
\usepackage{graphicx}
\usepackage{cite}
\usepackage{enumerate}
\usepackage{array}
\usepackage{float}
\usepackage{booktabs}
\usepackage{caption}
\usepackage{subcaption}
\usepackage{hyperref}
\usepackage{upgreek}
\usepackage{mathrsfs}
\hypersetup{
    pdftitle={Iterative Methods for Solving Fourth and Sixth Order Time-fractional Cahn-Hillard Equation},    % title
    pdfauthor={Lanre Akinyemi},     % author
    pdfcreator={pdftex},   % creator of the document
    pdfproducer={Texmaker}, % producer of the document
    pdfkeywords={Cahn-Hillard equation, Fractional derivatives, New iterative method, q-Homotopy analysis method}, % list of keywords
    colorlinks=true,       % false: boxed links; true: colored links
    linkcolor=red,          % color of internal links
    citecolor=blue,        % color of links to bibliography
    filecolor=magenta,      % color of file links
    urlcolor=blue           % color of external links
}
\DeclareMathOperator{\sech}{sech}
\newcommand{\capu}{\mathcal{D}}

\Crefname{figure}{Fig.}{Figs.}
\Crefname{table}{Tab.}{Tabs.}
\makeatletter
\let\oldtheequation\theequation
\renewcommand\tagform@[1]{\maketag@@@{\ignorespaces#1\unskip\@@italiccorr}}
\renewcommand\theequation{(\oldtheequation)}
\makeatother
\newtheorem{theorem}{Theorem}[section]

\newtheorem{definition}[theorem]{Definition}

\newtheorem{remark}[theorem]{Remark}

\date{}
\begin{document}

\title{Iterative Methods for Solving Fourth and Sixth Order Time-fractional Cahn-Hillard Equation}
\author{Lanre Akinyemi $^{1}$, Olaniyi S. Iyiola $^{2},\,$ Udoh Akpan$^{1}\,$
\\
\\
\\
\\
$^{1}$\textit{Department of Mathematics, Ohio University,}\\
\textit{Athens, OH 45701, United States of America}\\
$^{2}$\textit{Department of Mathematics, Computer Science \& Information System,}\\ 
\textit{California University of Pennsylvania, PA 15419, United States of America}
\\
\\
\\
e-mail: akinyemi@ohio.edu;\, iyiola@calu.edu;\; ua780216@ohio.edu}

\maketitle
%%%%%%%%%%%%%%%%%%%%%%%%%%%%%%%%%%%%%%%%%%%%%%%%%%%%%%%%%%%%%              %%%%%%%%%%% ABSTRACT %%%%%%%%%%%%%%
%%%%%%%%%%%%%%%%%%%%%%%%%%%%%%%%%%%%%%%%%%%%%%%%%%%%%%%%%%%%%

\begin{abstract}
This paper presents analytical-approximate solutions of the time-fractional Cahn-Hilliard (TFCH) equations of fourth and sixth-order using the new iterative method (NIM) and q-homotopy analysis method (q-HAM). We obtained convergent series solutions using these iterative methods. The simplicity and accuracy of these methods in solving strongly nonlinear fractional differential equations is displayed through the examples provided. In the case where exact solution exists, error estimates are also investigated.\\

\end{abstract}
%%%%%%%%%%%%%%%%%%%%%%%%%%%%%%%%%%%%%%%%%%%%%%%%%%%%%%%%%%%%%
          %%%%%%%%% Introduction %%%%%%%%%%%%%%%
%%%%%%%%%%%%%%%%%%%%%%%%%%%%%%%%%%%%%%%%%%%%%%%%%%%%%%%%%%%%%
\section{Introduction}
The concept of fractional calculus such as fractional derivatives and fractional integral is not new. \textit{L' Hospital}, in 1695, wrote a letter to \textit{Leibniz}, saying ``How do we calculate $d^n y/dx^n$ when $n=1/2$?" that is; ``what happen if we consider $n$ to be a fraction?" \textit{Leibniz} replied to \textit{L' Hospital} question saying ``$d^{1/2}x$ equal to $x\sqrt{dx:x}$" In actual fact , the reply is an apparent paradox, from this apparent paradox, one day useful result might be drawn \cite{Leib, Leib1, fc}. Later, researchers discovered that fractional derivative has a wide applications in the various fields of natural sciences and engineering such as found in control theory of dynamical systems, signal and image processing, financial modeling, nanotechnology, viscoelasticity, random walk, anomalous transport and anomalous diffusion, are just a few (for more details, see \cite{Bal, ChenChen, Ucs, Hilf1, Main, Monj,VEVV, Y, ZW, ZPZ, PU}).\\

Linear and nonlinear partial fractional differential equations has play a vital role in  various fields of engineering and natural sciences. Among such PDEs, we have Cahn-Hillard equation named after Cahn and Hilliard in 1958 \cite{CahnE}. This equation plays an essential role in understanding a number of fascinating physical phenomena for instance, in spinodal decomposition, phase ordering dynamics, and also describes vital qualitative distinctive attribute of two-phase systems connected with phase separation processes (see \cite{Cahn, CahnE, EllZ0, NOV} for a detailed discussion). Because of its real world applications in these various fields mention above, researchers have investigated the mathematical and numerical solutions of this equation \cite{DahBen, DeM, EllZ1, Jun, MM, NOV, RyHoff, UK}.\\

Solving partial differential equations with fractional derivatives is often more difficult than the classical PDEs, for it operator is defined by integral. In the recent year, researchers have develop some approximate methods in solving both linear and nonlinear PDEs with fractional derivatives, such as adomian decomposition method, (ADM) \cite{DahBen, UK, adm2}, variational iteration method, (VIM) \cite{Das, vim}, generalized differential transform method, (GDTM) \cite{Atangana1, Atangana2}, permuturbation iteration transformation method, (PITM) \cite{khals}, homotopy-perturbation method, (HPM) \cite{hpm, UK}, and  residual power series method (RPSM) \cite{Alquran2, Xu}.\\

In this paper we consider the time-fractional Cahn-Hilliard equations of the forth and sixth order given, respectively, by

\begin{equation}\label{eq:1}
 \capu_t^{\alpha}u=\mu u_x+(-u_{xx}-u+u^3)_{xx},
\end{equation}

and

\begin{equation}\label{eq:1N}
\capu_t^{\alpha}u=\mu uu_{x}+(u_{xx}+u-u^3)_{xxxx},
\end{equation}

with the initial condition
\begin{equation}\label{eq:1bc}
u(x,0)=f(x).
\end{equation}
Here, the parameters $\alpha$ ($0<\alpha \leq 1$) stands for the order of the time-fractional derivatives, satisfying $\mu \geq 0$. Our aim is to obtain solutions in the form of recurrence relations, using the new iteration method (NIM), which is base on the decomposition the non-linearity term \cite{Daf, Daf1, Bha, Wen} and q-homotopy analysis method (q-HAM), a modification of the homotopy analysis method (HAM) \cite{Taw, IyizOla, Iyiz1, Iyiz3}. 
%%%%%%%%%%%%%%%%%%%%%%%%%%%%%%%%%%%%%%%%%%%%%%%%%%%
          %%%%%% Basic Definitions %%%%%%%
%%%%%%%%%%%%%%%%%%%%%%%%%%%%%%%%%%%%%%%%%%%%%%%%%%%
\begin{definition}
The Gamma function is defined as \cite{Miller}
\begin{equation}\label{2}
\Gamma{(z)}=\int_{0}^{\infty} t^{z-1} e^{-t}\ dt,
\end{equation}
where $\Re(z)>0$.
\end{definition}

\begin{definition}
The Riemann-Liouville fractional integral operator of order $\alpha$ ($\alpha \geq 0$) of a function $u(x,t) \in C_{\mu }$, $\mu \geq -1$, denoted by $J^{\alpha }u(x,t)$ (with respect to ``t") is defined as 
\cite{fc,kil},
\begin{equation} \label{4}
J^{\alpha}u(x,  t)=\frac{1}{\Gamma (\alpha )}\int_{0}^{t}(t-\mu
)^{\alpha -1}\ u(x, \mu)\ d\mu ,\text{ \ \ }\alpha, \ t>0,
\end{equation}

where $J^{0}u(x,t)=u(x,t)$. Then the following properties holds for function $u(x,t)$ as follows:

\begin{enumerate}
\item[a.] $\displaystyle J^{\alpha }J^{\beta }u(x,t)=J^{\alpha +\beta }u(x,t)$

\item[b.] $\displaystyle J^{\alpha }J^{\beta}u(x,t)=J^{\beta }J^{\alpha }u(x,t)$

\item[c.] $\displaystyle J^{\alpha }t^{\uptau}=\frac{\Gamma (\uptau +1)}{\Gamma (\uptau
+1+ \alpha)}t^{\uptau+\alpha}$.
\end{enumerate}
Here, $\alpha$, $\beta \ \geq 0$ and $\uptau >-1.$
\end{definition}

\begin{definition}
The (left sided) Caputo fractional derivative of a function $u(x,t)$ of order $\alpha$ (with respect to ``t"), denoted by $\capu^{\alpha }u(x,t)$, where $m-1<\alpha <m,$ and $u(x,t) \in C_{l}^{m}$, $m \in \mathbb{N}$ is defined as \cite{fc,kil}
\begin{eqnarray}
\capu^{\alpha }u(x,t)=\left \{
\begin{array}{lll}
u^{(m)}(x,t),\qquad \qquad \alpha=m,\\
\\
J^{m-\alpha }u^{(m)}(x,t), \qquad m-1<\alpha < m,
\label{lem2}
\end{array}
\right.
\end{eqnarray}
where
\begin{equation}\label{5}
J^{m-\alpha }u^{(m)}(x, t)=\frac{1}{\Gamma (m-\alpha )}
\int_{0}^{t}(t-\xi)^{m-\alpha -1}u^{(m)}(x, \xi)\ d\xi ,\text{ \ \ }
\alpha,\ t>0,
\end{equation}
satisfies the following defined properties:
\begin{enumerate}
\item[a.] $ \capu^{\alpha}[\sigma\ u(x,t)+\rho \ w(x,t)] =\sigma \capu^{\alpha }u(x,  t)+\rho \capu^{\alpha
}w(x,  t),$ \quad $\sigma,\ \rho \in \mathbb{R} $,

\item[b.] $ \capu^{\alpha }J^{\alpha }u(x,t)=u(x,t)$,

\item[c.] $ J^{\alpha }\capu^{\alpha }u(x,  t)=u(x,  t)-\sum_{j=0}^{m-1}u^{j}(x, 0)\
\frac{t^{j}}{j!}$.
\end{enumerate}
\end{definition}
%%%%%%%%%%%%%%%%%%%%%%%%%%%%%%%%%%%%%%%%%%%%%%%%%%%%%%%%%%%%%%%%%
\section{Analysis of Approximate Methods}
%%%%%%%%%%%%%%%%%%%%%%%%%%%%%%%%%%%%%%%%%%%%%%%%%%%%%%%%%%%%%%%%%
In this subsection, we give a brief description of the new iterative method (NIM) and the q-homotopy analysis method (q-HAM). 
%%%%%%%%%%%%%%%%%%%%%%%%%%%%%%%%%%%%%%%%%%%%%%%%%%%%%%%%%%%
\subsection{Fundamentals of the NIM}
Consider the following functional equation,
\begin{equation}\label{eq:func}
 u(x,t)=f(x,t)+\mathscr{L}(u(x,t))+\mathcal{N}(u(x,t)),
\end{equation}
where $\mathscr{L}$ and $\mathcal{N}$ are respectively the linear and nonlinear operator from a Banach space $\mathcal{B}$ to itself, $f(x,t)$ is the known function and $u(x,t)$ is the unknown function. Let

\begin{eqnarray}\label{eq:func1}
U_n(x,t)=\sum_{k=0}^{n}u_k(x, t).
\end{eqnarray}
We seek for a solution $u(x, t)$ of \autoref{eq:func} has a series form 
\begin{eqnarray}\label{eq:func12}
u(x, t)=\lim_{n\to\infty}U_n(x,t)=\sum_{k=0}^{\infty} u_k(x, t).
\end{eqnarray}
\autoref{eq:func1} converges absolutely and uniformly to a unique solution if the operators $\mathscr{L}$ and $\mathcal{N}$ are contractive \cite{Daf, Wen}. The decomposed nonlinear operator $\mathcal{N}$ is give as 

\begin{eqnarray}\notag
\mathcal{N}(u(x,t))&=&\mathcal{N}\left(\sum_{k=0}^{\infty}u_k(x,t)\right)\\
&=&\mathcal{N}(u_0(x,t))+\sum_{k=1}^{\infty}\left \{\mathcal{N}\left(\sum_{i=0}^{k}u_i(x,t)\right)-\mathcal{N}\left(\sum_{i=0}^{k-1}u_i(x,t)\right)\right\}. \label{eq:func2}
\end{eqnarray}
In the same manner, the linear operator $\mathscr{L}$ can be decomposed as
\begin{eqnarray}\notag
\mathscr{L}(u(x,t))&=&\mathscr{L}\left(\sum_{k=0}^{\infty}u_k(x,t)\right)\\
&=&\mathscr{L}(u_0(x,t))+\sum_{k=1}^{\infty}\left\{\mathscr{L}\left(\sum_{i=0}^{k}u_i(x,t)\right)-\mathscr{L}\left(\sum_{i=0}^{k-1}u_i(x,t)\right)\right\}. \label{eq:func3}
\end{eqnarray}
Since $\mathscr{L}$ is a linear operator, then from \autoref{eq:func3} we obtain
\begin{eqnarray}\notag
\sum_{k=1}^{\infty}\left\{\mathscr{L}\left(\sum_{i=0}^{k}u_i(x,t)\right)-\mathscr{L}\left(\sum_{i=0}^{k-1}u_i(x,t)\right)\right\}&=&\mathscr{L}(u_1(x,t))+\mathscr{L}(u_2(x,t))+\mathscr{L}(u_3(x,t))+\cdots \\
&=&\sum_{k=1}^{\infty}\mathscr{L}(u_k(x,t)).\label{eq:func4}
\end{eqnarray}
Considering Eqs.~\ref{eq:func1}-\ref{eq:func4}, from \autoref{eq:func} we have
\begin{equation}\label{eq:func5}
\sum_{k=0}^{\infty}u_k(x,t)=f(x,t)+\sum_{k=0}^{\infty}\mathscr{L}(u_k(x,t))+\mathcal{N}(u_0(x,t))+\sum_{k=1}^{\infty}\left\{\mathcal{N}\left(\sum_{j=0}^{k}u_j(x,t)\right)-\mathcal{N}\left(\sum_{j=0}^{k-1}u_j(x,t)\right)\right\}.
\end{equation}
Then, from \autoref{eq:func5}, we define the iterations
\begin{eqnarray}\notag
u_0&=&f(x),\\ \notag
u_1&=&\mathscr{L}(u_0)+\mathcal{N}(u_0),\\ \notag
u_2&=&\mathscr{L}(u_1)+\big(\mathcal{N}(u_0+u_1)-\mathcal{N}(u_0)\big),\\ \notag
u_3&=&\mathscr{L}(u_2)+\big(\mathcal{N}(u_0+u_1+u_2)-\mathcal{N}(u_0+u_1)\big)\\ \notag
&\vdots&\\
u_{m+1}&=&\mathscr{L}(u_m)+\left\{\mathcal{N}\left(\sum_{j=0}^{m}u_k(x,t)\right)-\mathcal{N}\left(\sum_{j=0}^{m-1}u_k(x,t)\right)\right\}, \quad m=1,2,3,\cdots.\label{eq:func6}
\end{eqnarray}
%%%%%%%%%%%%%%%%%%%%%%%%%%%%%%%%%%%%%%%%%%%%%%%%
\subsection{Idea of the q-HAM}
Consider the nonlinear differential equation
\begin{equation}\label{eq:gen-eqn}
 \mathcal{N}\left(\capu_t^\alpha u(x, t)\right)-f(x, t)=0,
\end{equation}
where $f(x,t)$ and $u(x,t)$ are respectively the known and unknown functions, $\mathcal{N}$ is the nonlinear operator, and $\capu_t^\alpha$ is the conformable fractional derivative with respect to ``t". In order to generalize the concept of homotopy method, we construct the zeroth-order deformation equation given as
\begin{eqnarray}\label{eq:1a}
 (1-qn)\mathscr{L}\Big(\Phi(x,t;q)-u_0(x, t)\Big)=hq\mathcal{H}(x, t)\Big(\mathcal{N}[\capu_t^\alpha \Phi(x,t;q)]-f(x, t)\Big), \quad n\geq{1},
\end{eqnarray}
where $q$ ($0\leq q \leq\frac{1}{n}$) is the embedded parameter, $h$ ($h\neq0$) an auxiliary parameter, $\mathscr{L}$ is the auxiliary linear operator, and a non-zero auxiliary function denoted by $\mathcal{H}(x, t).$ For $q=0,\, \frac{1}{n}$ respectively, we obtain from \autoref{eq:1a} the following
\begin{equation}\label{eq:2a}
 \Phi\Big(x,t;0\Big)=u_0(x, t), \quad \quad \Phi\Big(x, t;\frac{1}{n}\Big)=u(x, t).
\end{equation}

When $q$ rises from $0$ to $\frac{1}{n}$, the solution $\Phi(x,t;q)$ ranges from the initial guess $u_0(x, t)$ to the solution $u(x, t)$. If $u_0(x, t)$, $h$, $\mathscr{L}$, and $\mathcal{H}(x, t)$ are chosen appropriately, then the solution $\Phi(x,t;q)$ of \autoref{eq:1a} is valid as long as $0\leq q \leq\frac{1}{n}$. Hence, we obtain the Taylor series expansion for $\Phi(x,t;q)$ as
\begin{equation}\label{eq:3a}
 \Phi(x,t;q)=u_0(x, t)+\sum_{r=1}^\infty u_r(x, t)q^r,
\end{equation}
where
\begin{equation}\label{eq:4a}
 u_r(x, t)=\left.\frac{1}{r!}\frac{\partial^r\Phi(x,t;q)}{\partial q^r}\right|_{q=0}.
\end{equation}

If we choose $u_0$, $h,$ $\mathscr{L}$, and $\mathcal{H}(x, t)$ properly so that the \autoref{eq:3a} converges at $q=\frac{1}{n}$, then from \autoref{eq:2a} we obtain
\begin{equation}
 u(x, t)=u_0(x, t)+\sum_{r=1}^\infty u_r(x, t)\left(\frac{1}{n}\right)^r.
\label{eq:5a}
\end{equation}

We define the vector $\vec{u}_r$ as follows:
\begin{equation}
\vec{u}_r=\left\{u_0(x, t),u_1(x, t),\cdots,u_r(x, t)\right\}.
\end{equation}

By differentiating \autoref{eq:1a} $r$-times (with respect to $``q"$), substituting $q=0,$ and then divide it by $r!$, we obtained what is known as the $r^{th}$-order deformation equation \cite{Liao, Iyiz3},
\begin{equation}
\mathscr{L}\left(u_r(x,t)-\Psi^*_r u_{r-1}(x, t)\right)=h\mathcal{H}(x, t)\mathcal{R}_r\left(\vec{u}_{r-1}\right).
\label{eq:6a}
\end{equation}
subject to the initial conditions\\
\begin{equation}
u_r^{(j)}(0, x)=0,\,\,\,\,\,\,\,\,\,j=0,1,2,3,... ,r-1.
\label{eq:7a}
\end{equation}
where
\begin{equation}\label{eq:8a}
\mathcal{R}_r\left(\vec{u}_{r-1}\right)=\left.\frac{1}{(r-1)!}\frac{\partial^{r-1}\left(\mathcal{N}[\capu_t^\alpha\Phi(x,t;q)]-f(x, t)\right)}{\partial q^{r-1}}\right|_{q=0}.
\end{equation}
Therefore, we have the solution to the system for $r\geq 1$ as
\begin{equation}\label{eq:6ss}
u_r(x, t)=\Psi^*_r u_{r-1}(x, t)+h J^\alpha\left(\mathcal{R}_r\left(\vec{u}_{r-1}(x, t)\right)\right),
\end{equation}
with
\begin{eqnarray}
\Psi^*_r=\left \{
\begin{array}{lll}
0\,\,\,\,\,\,\,\,\,\,\,\,\,r\leqslant 1,\\
\\
n\,\,\,\,\,\,\,\,\,\,\,\,\,otherwise.
\label{eq:9a}
\end{array}
\right.
\end{eqnarray}
The series solutions by q-HAM are
\begin{equation}\label{eq:sol1}
u(x,t;n;h)\cong U_r(x,t,n,h)=u_0(x, t)+\sum_{j=1}^r u_j(x,t;n;h)\left(\frac{1}{n}\right)^j,
\end{equation}
which gives the appropriate solutions in terms of convergence parameters $n$ and $h$.\\
%%%%%%%%%%%%%%%%%%%%%%%%%%%%%%%%%%%%%%%%%%%%%%%%%%%%%%%%%%%%%%%%%
        %%%%%%%%%%%%%%%% SOLUTIONS %%%%%%%%%%%%%%%%%%%%%%%
%%%%%%%%%%%%%%%%%%%%%%%%%%%%%%%%%%%%%%%%%%%%%%%%%%%%%%%%%%%%%%%%%
\section{Solutions of Fourth-Order TFCH Equation}\label{sec3}
In this section, we presents the application of the above mentioned methods to obtain approximate solution of the fourth-order TFCH \autoref{eq:1} subject to different initial conditions.
%%%%%%%%%%%%%%%%%%%%%%%%%%%%%%%%%%%%%%%%%%%%%%%%%%%%%%%%%%%%%%%%%
            %%%%%%%%%%%%%%% Case I %%%%%%%%%%%%%%%%%
%%%%%%%%%%%%%%%%%%%%%%%%%%%%%%%%%%%%%%%%%%%%%%%%%%%%%%%%%%%%%%%%%
\subsection{} 
\textbf{Case I:}
Consider the TFCH \autoref{eq:1}, 
\begin{eqnarray}\label{eq:7}
\capu_t^{\alpha}u=\mu u_x+6uu^2_x+3u^2u_{xx}-u_{xx}-u_{xxxx},\quad 0<\alpha \leq 1, \end{eqnarray}
with the initial condition
\begin{eqnarray}\label{eq:7IC}
u(x,0)=f(x)=\tanh{\left(\frac{x}{\sqrt{2}}\right)}.
\end{eqnarray}
The exact solution of the IVP \ref{eq:7}-\ref{eq:7IC}, when $\alpha, \,\mu= 1,$ is
\begin{equation}\label{eq:exactsol}
u(x,t)=\tanh{\left(\frac{x+t}{\sqrt{2}}\right)}.
\end{equation}\\
%%%%%%%%%%%%%%%%%%%%%%%%%%%%%%%%%%%%%%%%%%%%%%%%%%%%%%%%%%%%%%%%%
        %%%%%%%%%%%%%%%% NIM %%%%%%%%%%%%%%%%%%%
%%%%%%%%%%%%%%%%%%%%%%%%%%%%%%%%%%%%%%%%%%%%%%%%%%%%%%%%%%%%%%%%%
\textbf{NIM solution:}\\
Applying $J^\alpha$ to both sides of \autoref{eq:7}, then the IVP \ref{eq:7}-\ref{eq:7IC} is equivalent to the integral equation:
\begin{eqnarray*}
u(x, t)&=&f(x, t)+\mathscr{L}(u(x, t))+\mathcal{N}(u(x, t)),
\end{eqnarray*}
where, 
\begin{eqnarray*}
u_0&=&f(x)=\tanh{\left(\frac{x}{\sqrt{2}}\right)},\\
\mathscr{L}(u)&=&J^\alpha\big(\mu u_x -u_{xx}-u_{xxxx}\big),\\
\mathcal{N}(u)&=&J^\alpha\left(6 u u^2_{x}+3u^2u_{xx}\right).
\end{eqnarray*}
We now obtain components of the series solution using NIM recurrent relation in \autoref{eq:func6} successively as follows:
\begin{eqnarray*}
u_1(x,t)&=& \mathscr{L}(u_0)+\mathcal{N}(u_0)\\
&=&J^\alpha\big(\mu u_{0x}-u_{0xx}-u_{0xxxx}\big)\\
&&+\ J^\alpha\big(6u_0u^2_{0x}+3u^2_{0}u_{0xx}\big)\\
&=& \frac{\mu \sech^2{\left(\frac{x}{\sqrt{2}}\right)}}{\sqrt{2}\Gamma(\alpha+1)}t^{\alpha},
\end{eqnarray*}

\begin{eqnarray*}
u_2(x,t)&=& \mathscr{L}(u_1)+\big\{\mathcal{N}(u_0+u_1)-\mathcal{N}(u_0)\big\}\\
&=& J^\alpha\left(\mu u_{1x}-u_{1xx}-u_{1xxxx}+12u_0u_{0x}u_{1x}+6u_{1}u^2_{0x}+6u_{0}u_{1}u_{0xx}+3u^2_{0}u_{1xx}\right)\\
&&+J^\alpha\left(6u_0u^2_{1x}+12u_{1}u_{0x}u_{1x}+3u^2_{1}u_{0xx}+6u_{0}u_{1}u_{1xx}+6u_{1}u^2_{1x}+3u^2_{1}u_{1xx}\right)\\
&=& -\frac{\mu^2 \tanh{\left(\frac{x}{\sqrt{2}}\right)}\sech^2{\left(\frac{x}{\sqrt{2}}\right)}}{\Gamma(2\alpha+1)}t^{2\alpha}+\frac{3 \mu^2\Gamma(2\alpha+1)\left(4\cosh(\sqrt{2}x)-11\right)\tanh\left(\frac{x}{\sqrt{2}}\right)\sech^6{\left(\frac{x}{\sqrt{2}}\right)}}{2\Gamma(\alpha+1)^2\Gamma(3\alpha+1)}t^{3\alpha}\\
&&+\frac{3\mu^3\Gamma(3\alpha+1)\left(3\cosh(\sqrt{2}x)-4\right)\sech^8{\left(\frac{x}{\sqrt{2}}\right)}}{2\sqrt{2}\Gamma(\alpha+1)^3\Gamma(4\alpha+1)}t^{4\alpha}. 
\end{eqnarray*}
Using the same procedure, expression for $u_m(x,t),$ $m=3,4,5,...$ can be obtained. The expression of the series solutions given by this new iteration method (NIM) can be written in the form

\begin{eqnarray}\notag
u(x,t)&\cong & \sum_{i=0}^{2} u_i(x,t)=u_0(x,t)+u_1(x,t)+u_2(x,t)\\ \notag
&=& \tanh{\left(\frac{x}{\sqrt{2}}\right)}+\frac{\mu \sech^2{\left(\frac{x}{\sqrt{2}}\right)}}{\sqrt{2}\Gamma(\alpha+1)}t^\alpha-\frac{\mu^2 \tanh{\left(\frac{x}{\sqrt{2}}\right)}\sech^2{\left(\frac{x}{\sqrt{2}}\right)}}{\Gamma(2\alpha+1)}t^{2\alpha}\\ \notag
&&+\frac{3\mu^2\Gamma(2\alpha+1)\left(4\cosh(\sqrt{2}x)-11\right)\tanh\left(\frac{x}{\sqrt{2}}\right)\sech^6{\left(\frac{x}{\sqrt{2}}\right)}}{2\Gamma(\alpha+1)^2\Gamma(3\alpha+1)}t^{3\alpha}\\ \label{eq:solutnimm}
&&+\frac{3\mu^3\Gamma(3\alpha+1)\left(3\cosh(\sqrt{2}x)-4\right)\sech^8{\left(\frac{x}{\sqrt{2}}\right)}}{2\sqrt{2}\Gamma(\alpha+1)^3\Gamma(4\alpha+1)}t^{4\alpha}. 
\end{eqnarray}
Thus, \autoref{eq:solutnimm} gives an approximate solution to the IVP \ref{eq:7}-\ref{eq:7IC}.\\
%%%%%%%%%%%%%%%%%%%%%%%%%%%%%%%%%%%%%%%%%%%%%%%%%%%%%%%%%%%%%%%%%
        %%%%%%%%%%%%%%%% q-HAM %%%%%%%%%%%%%%%%%%%
%%%%%%%%%%%%%%%%%%%%%%%%%%%%%%%%%%%%%%%%%%%%%%%%%%%%%%%%%%%%%%%%%
\newline
\textbf{q-HAM solution:}\\
To apply the q-HAM, we rewrite TFCH \autoref{eq:1} as
\begin{eqnarray}\label{eq:7q}
\capu_t^{\alpha}u-\mu u_x-6uu^2_x-3u^2u_{xx}+u_{xx}+u_{xxxx}=0,\quad 0<\alpha \leq 1. \end{eqnarray}
 Applying q-HAM to \autoref{eq:7q}, we obtain from \autoref{eq:8a},
\begin{eqnarray}\notag
\mathcal{R}_m\left(\vec{u}_{m-1}\right)&=&\capu_t^\alpha u_{m-1}-\mu u_{(m-1)x}-6\sum_{k=0}^{m-1}\sum_{j=0}^{k}u_j u_{(k-j)x}u_{(m-1-k)x}\\
&& - \ 3\sum_{k=0}^{m-1}\sum_{j=0}^{k}u_j u_{(k-j)}u_{(m-1-k)xx} + u_{(m-1)xx} + u_{(m-1)xxxx}.
\label{eq:6s}
\end{eqnarray}
We can now obtain the components of the solution using q-HAM recurrent relation in \autoref{eq:6ss}, using Eqs.~\ref{eq:9a} and \ref{eq:6s} successively as follows:
\begin{eqnarray*}
u_1(x,t)&=&\Psi^*_1 u_0(x,t)+h J^\alpha\left(\mathcal{R}_1\left(\vec{u}_0(x,t)\right)\right)\\
&=& h J^\alpha\left(\mathcal{D}_t^\alpha u_0 - \mu u_{0x}-6u_0u_{0x}u_{0x}\right)\\
&&+\ h J^\alpha\left(-3u_{0}u_{0}u_{0xx}+u_{0xx}+u_{0xxxx}\right)\\
&=& -\frac{h\mu \sech^2{\left(\frac{x}{\sqrt{2}}\right)}}{\sqrt{2}\Gamma(\alpha+1)}t^\alpha ,
\end{eqnarray*}

\begin{eqnarray*}
u_2(x,t)&=&\Psi^*_2 u_1(x,t)+h J^\alpha\left(\mathcal{R}_2\left(\vec{u}_1(x,t)\right)\right)\\
&=&hu_1+ h J^\alpha\left(\mathcal{D}_t^\alpha u_1 - \mu u_{1x}-6u_0u_{0x}u_{1x}-6u_0u_{1x}u_{0x}-6u_1u_{0x}u_{0x}\right)\\
&&+ \ h J^\alpha\left(-3u_{0}u_{0}u_{1xx}-3u_{0}u_{1}u_{0xx}-3u_{1}u_{0}u_{0xx}+u_{1xx}+u_{1xxxx}\right)\\
&=& (n+h)\,u_1-\frac{h^2 \mu^2 \tanh{\left(\frac{x}{\sqrt{2}}\right)}\sech^2{\left(\frac{x}{\sqrt{2}}\right)}}{\Gamma(2\alpha+1)}t^{2\alpha},
\end{eqnarray*}

\begin{eqnarray*}
u_3(x,t)&=&\Psi^*_3 u_2(x,t)+h J^\alpha\left(\mathcal{R}_3\left(\vec{u}_2(x,t)\right)\right)\\
&=&hu_2+ h J^\alpha\left(\mathcal{D}_t^\alpha u_2 - \mu u_{2x}-6u_0u_{0x}u_{2x}-6u_0u_{1x}u_{1x}-6u_1u_{0x}u_{1x}\right)\\
&&+\ h J^\alpha\left(-6u_0u_{2x}u_{0x}-6u_1u_{1x}u_{0x}-6u_2u_{0x}u_{0x}-3u_{0}u_{0}u_{2xx}-3u_{0}u_{1}u_{1xx}\right)\\
&&+ \ h J^\alpha\left(-3u_{1}u_{0}u_{1xx}-3u_{0}u_{2}u_{0xx}-3u_{1}u_{1}u_{0xx}-3u_{2}u_{0}u_{0xx}+u_{2xx}+u_{2xxxx}\right)\\
&=& (n+h)\,u_2-\frac{h^2\mu^2(n+h) \tanh{\left(\frac{x}{\sqrt{2}}\right)}\sech^2{\left(\frac{x}{\sqrt{2}}\right)}}{\Gamma(2\alpha+1)}t^{2\alpha}\\
&&-\frac{\sqrt{2}h^3 \mu^3\Big(\cosh(\sqrt{2}x)-2\Big)\sech^4{\left(\frac{x}{\sqrt{2}}\right)}}{2\Gamma(3\alpha+1)}t^{3\alpha}\\
&&+\frac{6h^3 \mu^2\Big(4\cosh(\sqrt{2}x)-11\Big)\tanh\left(\frac{x}{\sqrt{2}}\right)\sech^6{\left(\frac{x}{\sqrt{2}}\right)}}{2\Gamma(3\alpha+1)}t^{3\alpha}\\%%%
&&- \frac{3h^3 \mu^2\Gamma(2\alpha+1)\left(4\cosh(\sqrt{2}x)-11\right)\tanh\left(\frac{x}{\sqrt{2}}\right)\sech^6{\left(\frac{x}{\sqrt{2}}\right)}}{2\Gamma(\alpha+1)^2\Gamma(3\alpha+1)}t^{3\alpha}.
\end{eqnarray*}

Using the same procedure, expression for $u_m(x,t),$ $m=4,5,6,...$ can be obtained. The expression of the series solutions given by q-HAM can be written in the form

\begin{eqnarray}\notag
u(x,t;n;h)&\cong&u_0(x,t)+\sum_{i=1}^3 u_i(x,t;n;h)\left(\frac{1}{n}\right)^i\\ \notag
&=& \tanh{\left(\frac{x}{\sqrt{2}}\right)}-\frac{\mu h(3n^2+3nh+h^2) \sech^2{\left(\frac{x}{\sqrt{2}}\right)}}{\sqrt{2}n^3\Gamma(\alpha+1)}t^\alpha-\frac{h^2 \mu^2 (3n+2h) \tanh{\left(\frac{x}{\sqrt{2}}\right)}\sech^2{\left(\frac{x}{\sqrt{2}}\right)}}{n^3\Gamma(2\alpha+1)}t^{2\alpha}\\ \notag
&&-\frac{h^3 \mu^2\sech^4{\left(\frac{x}{\sqrt{2}}\right)}\left(\sqrt{2}\mu\left(\cosh(\sqrt{2}x)-2\right)-6(4\cosh(\sqrt{2}x-11)\tanh\left(\frac{x}{\sqrt{2}}\right)\sech^2{\left(\frac{x}{\sqrt{2}}\right)}\right)}{2n^3\Gamma(3\alpha+1)}t^{3\alpha}\\ \label{eq:solut}
&&-\frac{3h^3 \mu^2\Gamma(2\alpha+1)\left(4\cosh(\sqrt{2}x)-11\right)\tanh\left(\frac{x}{\sqrt{2}}\right)\sech^6{\left(\frac{x}{\sqrt{2}}\right)}}{2n^3\Gamma(\alpha+1)^2\Gamma(3\alpha+1)}t^{3\alpha}.
\end{eqnarray}
Thus, \autoref{eq:solut} gives an approximate solutions to the IVP \ref{eq:7}-\ref{eq:7IC} in terms of convergence parameters $h$ and $n$. In the case when $n=\alpha=1,$ we choose $h=-1,$ and obtain from \autoref{eq:solut} the expression
\begin{eqnarray}\notag
u(x,t)&=&\tanh{\left(\frac{x}{\sqrt{2}}\right)}+\frac{\mu \sech^2{\left(\frac{x}{\sqrt{2}}\right)}}{\sqrt{2}}t-\frac{\mu^2 \tanh{\left(\frac{x}{\sqrt{2}}\right)}\sech^2{\left(\frac{x}{\sqrt{2}}\right)}}{2}t^{2} \\
&& + \frac{\mu^3 \left(\cosh(\sqrt{2}x)-2\right)\sech^4{\left(\frac{x}{\sqrt{2}}\right)}}{6\sqrt{2}}t^{3}+\cdots, \label{solution}
\end{eqnarray}
which can be expressed in the closed form of the exact solution
\begin{equation}\label{exacts}
u(x,t)=\tanh{\left(\frac{x+t}{\sqrt{2}}\right)}.
\end{equation}
\begin{remark}
This agrees with the solution obtained using ADM and HPM in \cite{UK}.
\end{remark}
%%%%%%%%%%%%%%%%%%%%%%%%%%%%%%%%%%%%%%%%%%%%%%%%%%%%%%%%%%%%%%%%%
         %%%%%%%%%%%%%%%%% cASE II %%%%%%%%%%%%%%%%%%%%
%%%%%%%%%%%%%%%%%%%%%%%%%%%%%%%%%%%%%%%%%%%%%%%%%%%%%%%%%%%%%%%%
\newpage
\textbf{Case II:}
Consider the TFCHE \autoref{eq:1}, 
\begin{eqnarray}\label{eq:10d}
\capu_t^{\alpha}u=\mu u_x+6uu^2_x+3u^2u_{xx}-u_{xx}-u_{xxxx}, \quad 0<\alpha \leq 1,
\end{eqnarray}
with the initial condition
\begin{eqnarray}\label{eq:10dIC}
u(x,0)=f(x)=e^{\lambda x}.
\end{eqnarray}
%%%%%%%%%%%%%%%%%%%%%%%%%%%%%%%%%%%%%%%%%%%%%%%%%%%%%%%%%%%%%%%%%
          %%%%%%%%%%%%%%%%%% NIM %%%%%%%%%%%%%%%%%%
%%%%%%%%%%%%%%%%%%%%%%%%%%%%%%%%%%%%%%%%%%%%%%%%%%%%%%%%%%%%%%%%%
\newline
\textbf{NIM solution:}\\
Applying $J^\alpha$ to both sides of \autoref{eq:10d}, then the IVP \ref{eq:10d}-\ref{eq:10dIC} is equivalent to the integral equation:
\begin{eqnarray*}
u(x, t)&=& f(x, t)+\mathscr{L}(u(x, t))+\mathcal{N}(u(x, t)),
\end{eqnarray*}
where, 
\begin{eqnarray*}
u_0&=&f(x)=e^{\lambda x},\\
\mathscr{L}(u)&=&J^\alpha\left(\mu u_x -u_{xx}-u_{xxxx}\right),\\
\mathcal{N}(u)&=&J^\alpha\left(6 u u^2_{x}+3u^2u_{xx}\right).
\end{eqnarray*}
We now obtain components of the series solution using NIM recurrent relation in \autoref{eq:func6} successively as follows:
\begin{eqnarray*}
u_1(x,t)&=& \mathscr{L}(u_0)+\mathcal{N}(u_0)\\
&=& J^\alpha\left(\mu u_{0x}-u_{0xx}-u_{0xxxx}\right)\\
&&+\ J^\alpha\left(6u_0u^2_{0x}+3u^2_{0}u_{0xx}\right)\\
&=& \frac{\lambda e^{\lambda x}(- \lambda^3+9 \lambda e^{2\lambda x}- \lambda+\mu) }{\Gamma(\alpha+1)}t^\alpha ,
\end{eqnarray*}

\begin{eqnarray*}
u_2(x,t)&=& \mathscr{L}(u_1)+\big\{\mathcal{N}(u_0+u_1)-\mathcal{N}(u_0)\big\}\\
&=& J^\alpha\left(\mu u_{1x}-u_{1xx}-u_{1xxxx}+12u_0u_{0x}u_{1x}+6u_{1}u^2_{0x}+6u_{0}u_{1}u_{0xx}+3u^2_{0}u_{1xx}\right)\\
&&+J^\alpha\left(6u_0u^2_{1x}+12u_{1}u_{0x}u_{1x}+3u^2_{1}u_{0xx}+6u_{0}u_{1}u_{1xx}+6u_{1}u^2_{1x}+3u^2_{1}u_{1xx}\right)\\
&=& \frac{\lambda^2e^{\lambda x}\left(-54\lambda e^{2\lambda x}(14\lambda^3+2\lambda-\mu)+(\lambda^3+\lambda-\mu)^2+675\lambda^2e^{4\lambda x}\right)}{\Gamma(2\alpha+1)}t^{2\alpha}\\
&&+\frac{27\lambda^4h^3 e^{3\lambda x}\Gamma(2\alpha+1)\left(50\lambda(\lambda^3+\lambda-\mu)e^{2\lambda x}-(\lambda^3+\lambda-\mu)^2-441\lambda^2e^{4\lambda x}\right)}{\Gamma(\alpha+1)^2\Gamma(3\alpha+1)}t^{3\alpha}\\
&&-\frac{9\lambda^5e^{3\lambda x}\Gamma(3\alpha+1)\left(\lambda^9+3\lambda^7+3\lambda^5+\lambda^3(3\mu^2+1)-3\lambda\mu(\lambda^5+2\lambda^3+\lambda-\mu)+1323(\lambda^5+\lambda^3-\lambda^2\mu)e^{4\lambda x}\right)}{\Gamma(\alpha+1)^3\Gamma(4\alpha+1)}t^{4\alpha}\\
&&+\frac{9\lambda^5e^{3\lambda x}\Gamma(3\alpha+1)\left(75(\lambda^7-2\lambda^4\mu-2\lambda^2\mu+\lambda\mu^2+2\lambda^5+\lambda^3)e^{2\lambda x}+6561\lambda^3e^{6\lambda x}-\mu^3 \right)}{\Gamma(\alpha+1)^3\Gamma(4\alpha+1)}t^{4\alpha}.
\end{eqnarray*}
Using the same procedure, expression for $u_m(x,t),$ $m=3,4,5,...$ can be obtained. The expression of the series solutions given by NIM can be written in the form

\begin{eqnarray}\notag
u(x,t)&\cong&\sum_{i=0}^{2} u_i(x,t)=u_0(x,t)+u_1(x,t)+u_2(x,t)\\ \notag
&=& e^{\lambda x}+\frac{\lambda e^{\lambda x}(- \lambda^3+9 \lambda e^{2\lambda x}- \lambda+\mu) }{\Gamma(\alpha+1)}t^\alpha\\ \notag
&&+ \frac{\lambda^2e^{\lambda x}\left(-54\lambda e^{2\lambda x}(14\lambda^3+2\lambda-\mu)+(\lambda^3+\lambda-\mu)^2+675\lambda^2e^{4\lambda x}\right)}{\Gamma(2\alpha+1)}t^{2\alpha}\\ \notag
&&+\frac{27\lambda^4h^3 e^{3\lambda x}\Gamma(2\alpha+1)\left(50\lambda(\lambda^3+\lambda-\mu)e^{2\lambda x}-(\lambda^3+\lambda-\mu)^2-441\lambda^2e^{4\lambda x}\right)}{\Gamma(\alpha+1)^2\Gamma(3\alpha+1)}t^{3\alpha}\\ \notag
&&-\frac{9\lambda^5e^{3\lambda x}\Gamma(3\alpha+1)\left(\lambda^9+3\lambda^7+3\lambda^5+\lambda^3(3\mu^2+1)-3\lambda\mu(\lambda^5+2\lambda^3+\lambda-\mu)+1323(\lambda^5+\lambda^3-\lambda^2\mu)e^{4\lambda x}\right)}{\Gamma(\alpha+1)^3\Gamma(4\alpha+1)}t^{4\alpha}\\ \label{eq:solutnim1}
&&+\frac{9\lambda^5e^{3\lambda x}\Gamma(3\alpha+1)\left(75(\lambda^7-2\lambda^4\mu-2\lambda^2\mu+\lambda\mu^2+2\lambda^5+\lambda^3)e^{2\lambda x}+6561\lambda^3e^{6\lambda x}-\mu^3 \right)}{\Gamma(\alpha+1)^3\Gamma(4\alpha+1)}t^{4\alpha}.
\end{eqnarray}
Thus, \autoref{eq:solutnim1} gives an approximate solution to the IVP \ref{eq:10d}-\ref{eq:10dIC}.\\
%%%%%%%%%%%%%%%%%%%%%%%%%%%%%%%%%%%%%%%%%%%%%%%%%%%%%%%%%%%%%%%%%
          %%%%%%%%%%%%%%%%%% q-HAM %%%%%%%%%%%%%%%%%%
%%%%%%%%%%%%%%%%%%%%%%%%%%%%%%%%%%%%%%%%%%%%%%%%%%%%%%%%%%%%%%%%%
\newline
\textbf{q-HAM solution:}\\
Consider \autoref{eq:7q}, we obtain the components of the solution using q-HAM recurrent relation in  \autoref{eq:6ss}, using Eqs.~\ref{eq:9a} and \ref{eq:6s} successively as follows:
\begin{eqnarray*}
u_1(x,t)&=&\Psi^*_1 u_0(x,t)+h J^\alpha\left(\mathcal{R}_1\left(\vec{u}_0(x,t)\right)\right)\\
&=& h J^\alpha\left(\mathcal{D}_t^\alpha u_0 - \mu u_{0x}-6u_0u_{0x}u_{0x}\right)\\
&&+\ h J^\alpha\left(-3u_{0}u_{0}u_{0xx}+u_{0xx}+u_{0xxxx}\right)\\
&=& -\frac{ \lambda h e^{\lambda x}(- \lambda^3+9 \lambda e^{2\lambda x}- \lambda+\mu) }{\Gamma(\alpha+1)}t^\alpha
\end{eqnarray*}

\begin{eqnarray*}
u_2(x,t)&=&\Psi^*_2 u_1(x,t)+h J^\alpha\left(\mathcal{R}_2\left(\vec{u}_1(x,t)\right)\right)\\
&=&hu_1+ h J^\alpha\left(\mathcal{D}_t^\alpha u_1 - \mu u_{1x}-6u_0u_{0x}u_{1x}-6u_0u_{1x}u_{0x}-6u_1u_{0x}u_{0x}\right)\\
&&+ \ h J^\alpha\left(-3u_{0}u_{0}u_{1xx}-3u_{0}u_{1}u_{0xx}-3u_{1}u_{0}u_{0xx}+u_{1xx}+u_{1xxxx}\right)\\
&=&(n+h)\,u_1+\frac{\lambda^2h^2e^{\lambda x}\left(-54\lambda e^{2\lambda x}(14\lambda^3+2\lambda-\mu)+(\lambda^3+\lambda-\mu)^2+675\lambda^2e^{4\lambda x}\right)}{\Gamma(2\alpha+1)}t^{2\alpha}
\end{eqnarray*}

\begin{eqnarray*}
u_3(x,t)&=&\Psi^*_3 u_2(x,t)+h J^\alpha\left(\mathcal{R}_3\left(\vec{u}_2(x,t)\right)\right)\\
&=&hu_2+ h J^\alpha\left(\mathcal{D}_t^\alpha u_2 - \mu u_{2x}-6u_0u_{0x}u_{2x}-6u_0u_{1x}u_{1x}-6u_1u_{0x}u_{1x}-6u_0u_{2x}u_{0x}-6u_1u_{1x}u_{0x}-6u_2u_{0x}u_{0x}\right)\\
&&+\ h J^\alpha\left(-3u_{0}u_{0}u_{2xx}-3u_{0}u_{1}u_{1xx}-3u_{1}u_{0}u_{1xx}-3u_{0}u_{2}u_{0xx}-3u_{1}u_{1}u_{0xx}-3u_{2}u_{0}u_{0xx}+u_{2xx}+u_{2xxxx}\right)\\
&=& (n+h)\,u_2+\frac{\lambda^2h^2(n+h)e^{\lambda x}\left(-54\lambda e^{2\lambda x}(14\lambda^3+2\lambda-\mu)+(\lambda^3+\lambda-\mu)^2+675\lambda^2e^{4\lambda x}\right)}{\Gamma(2\alpha+1)}t^{2\alpha}\\
&&+\frac{\lambda^3h^3e^{\lambda x}\left((\lambda^3+\lambda-\mu)^3-99225\lambda^3e^{6\lambda x}+675\lambda^3(709\lambda^2+
37)e^{4\lambda x}-27\lambda^3(2269\lambda^4+578\lambda^2+37
)e^{2\lambda x}
\right)}{\Gamma(3\alpha+1)}t^{3\alpha}\\
&&-\frac{27\mu\lambda^4h^3e^{3\lambda x} \left(-248\lambda^3+\lambda(275e^{2\lambda x}-32)+7\mu\right)}{\Gamma(3\alpha+1)}t^{3\alpha}\\
&&+\frac{27\lambda^4h^3 e^{3\lambda x}\Gamma(2\alpha+1)\left(50\lambda(\lambda^3+\lambda-\mu)e^{2\lambda x}-(\lambda^3+\lambda-\mu)^2-441\lambda^2e^{4\lambda x}\right)}{\Gamma(\alpha+1)^2\Gamma(3\alpha+1)}t^{3\alpha}
\end{eqnarray*}
Using the same procedure, expression for $u_m(x,t),$ $m=4,5,6,...$ can be obtained. The expression of the series solutions given by q-HAM can be written in the form

\begin{eqnarray}\notag 
u(x,t;n;h)&\cong& u_0(x,t)+\sum_{i=1}^3 u_i(x,t;n;h)\left(\frac{1}{n}\right)^i \\ \notag
&=&e^{\lambda x}-\frac{ \lambda{e^{\lambda x}}h(3n^2+3nh+h^2)(- \lambda^3+9 \lambda e^{2\lambda x}- \lambda+\mu)}{n^3\Gamma(\alpha+1)}t^\alpha \\ \notag
&&+\frac{\lambda^2h^2(3n+2h)e^{\lambda x}\left(-54\lambda e^{2\lambda x}(14\lambda^3+2\lambda-\mu)+(\lambda^3+\lambda-\mu)^2+675\lambda^2e^{4\lambda x}\right)}{n^3\Gamma(2\alpha+1)}t^{2\alpha}\\ \notag
&&+\frac{\lambda^3h^3e^{\lambda x}\left((\lambda^3+\lambda-\mu)^3-99225\lambda^3e^{6\lambda x}+675\lambda^3(709\lambda^2+
37)e^{4\lambda x}-27\lambda^3(2269\lambda^4+578\lambda^2+37
)e^{2\lambda x}\right)}{n^3\Gamma(3\alpha+1)}t^{3\alpha}\\ \notag
&&-\frac{27\mu\lambda^4h^3e^{3\lambda x} \left(-248\lambda^3+\lambda(275e^{2\lambda x}-32)+7\mu\right)}{n^3\Gamma(3\alpha+1)}t^{3\alpha}\\ \label{eq:solut5}
&&+\frac{27\lambda^4h^3 e^{3\lambda x}\left(50\lambda(\lambda^3+\lambda-\mu)e^{2\lambda x}-(\lambda^3+\lambda-\mu)^2-441\lambda^2e^{4\lambda x}\right)}{n^3\Gamma(\alpha+1)^2\Gamma(3\alpha+1)}t^{3\alpha}.
\end{eqnarray}
Thus, \autoref{eq:solut5} give an approximate solution to the IVP \ref{eq:10d}-\ref{eq:10dIC} in terms of convergence parameters $h$ and $n$.
%%%%%%%%%%%%%%%%%%%%%%%%%%%%%%%%%%%%%%%%%%%%%%%%%%%%%
      %%%%%%%%% NUMERICAL SOLUTION %%%%%%%%%%%
%%%%%%%%%%%%%%%%%%%%%%%%%%%%%%%%%%%%%%%%%%%%%%%%%%%%%
\subsection{Numerical Results for TFCH Equation of Fourth-Order}
\begin{figure}[H]
    \centering
    \begin{subfigure}[b]{0.32\textwidth}
        \includegraphics[width=\textwidth]{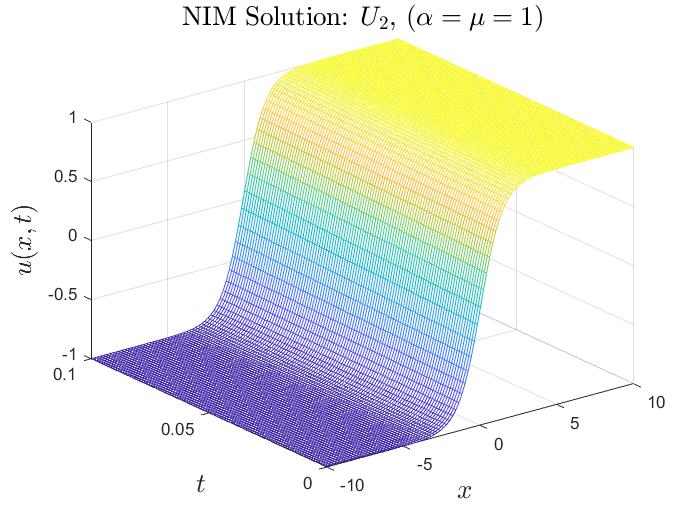}
        %\caption{}
        %\label{Fig1u}
    \end{subfigure}
    ~ %add desired spacing between images, e. g. ~, \quad, \qquad, \hfill etc. 
    \begin{subfigure}[b]{0.32\textwidth}
        \includegraphics[width=\textwidth]{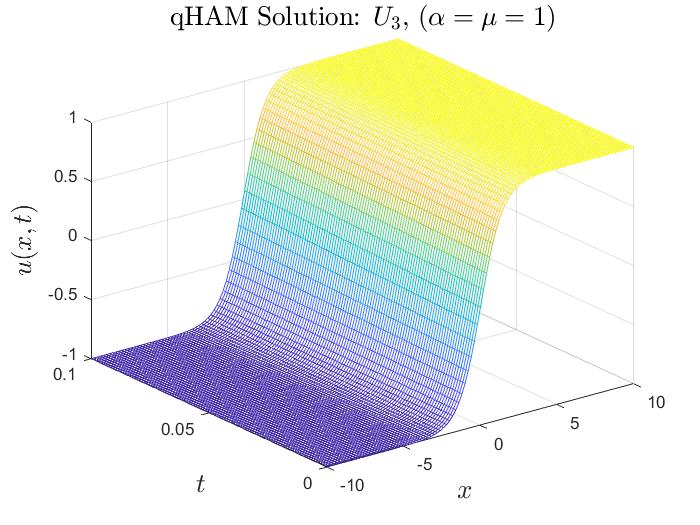}
        %\caption{}
        %\label{fi2}
    \end{subfigure}
    ~ %add desired spacing between images, e. g. ~, \quad, \qquad, \hfill etc. 
    \begin{subfigure}[b]{0.32\textwidth}
        \includegraphics[width=\textwidth]{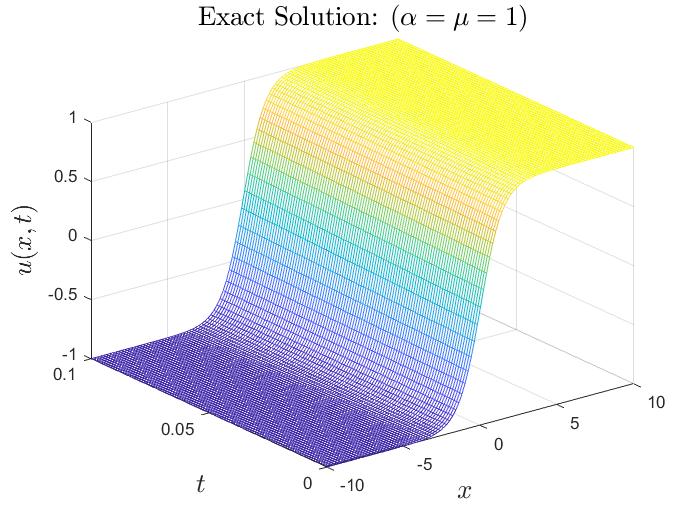}
       % \caption{}
       % \label{fi3}
    \end{subfigure}
\caption{Case I: Comparison among NIM, q-HAM and Exact Solution for $n=1,$ and $h=-1.$}\label{fig:Fig1}
\end{figure}

\begin{figure}[H]
\centering
    \begin{subfigure}[b]{0.40\textwidth}
        \includegraphics[width=\textwidth]{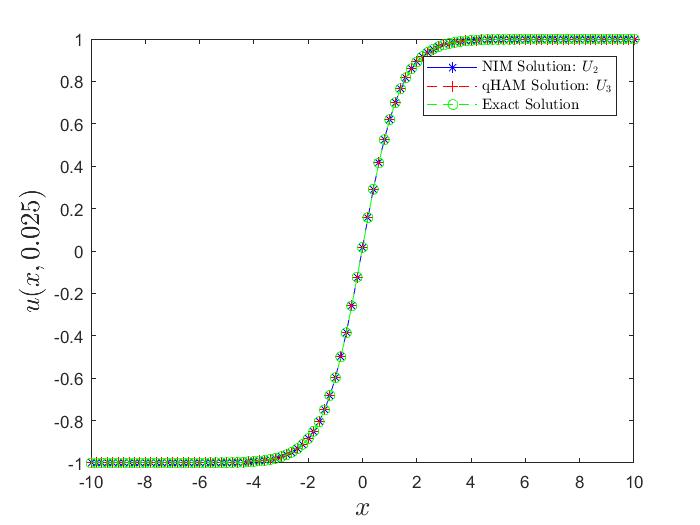}
        %\caption{}
    \end{subfigure}
\caption{Case I: Comparison among NIM, q-HAM and Exact Solution in 2D for $n=\alpha=\mu=1$ and $h=-1.$}\label{fig:Fig2}
\end{figure}

\begin{figure}[H]
    \centering
    \begin{subfigure}[b]{0.40\textwidth}
        \includegraphics[width=\textwidth]{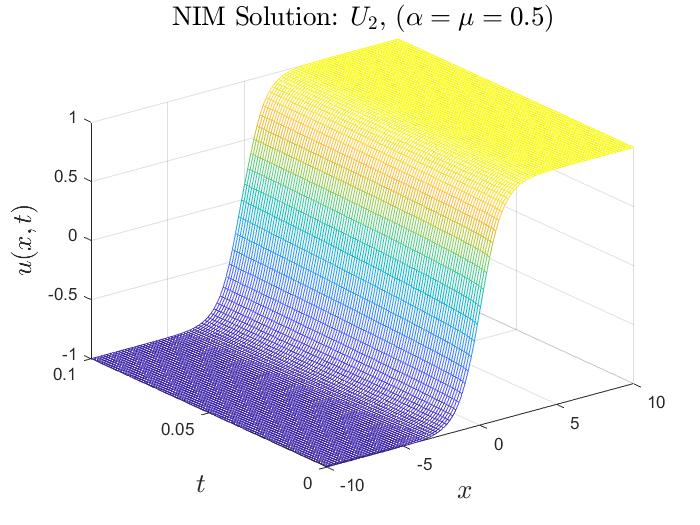}
       % \caption{}
        %\label{Fig1u}
    \end{subfigure}
    ~ %add desired spacing between images, e. g. ~, \quad, \qquad, \hfill etc. 
    \begin{subfigure}[b]{0.40\textwidth}
        \includegraphics[width=\textwidth]{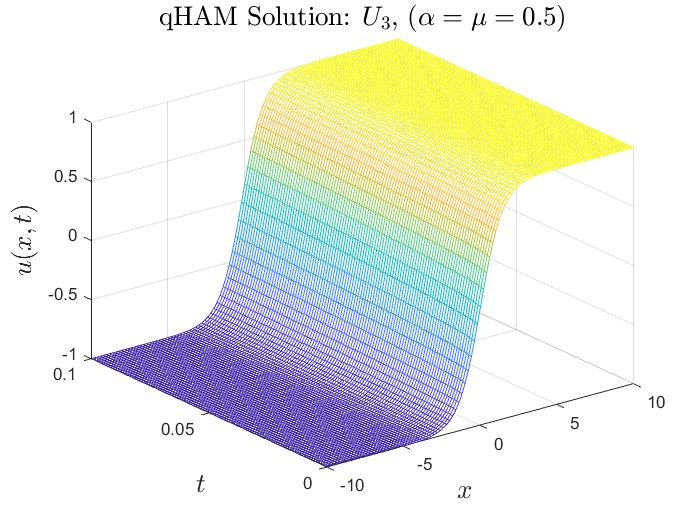}
        %\caption{}
        %\label{fi2}
    \end{subfigure}
\caption{Case I: Comparison between NIM and q-HAM Solution for $n=1$ and $h=-1.$}\label{fig:Fig3}
\end{figure}

\begin{figure}[H]
    \centering
    \begin{subfigure}[b]{0.40\textwidth}
        \includegraphics[width=\textwidth]{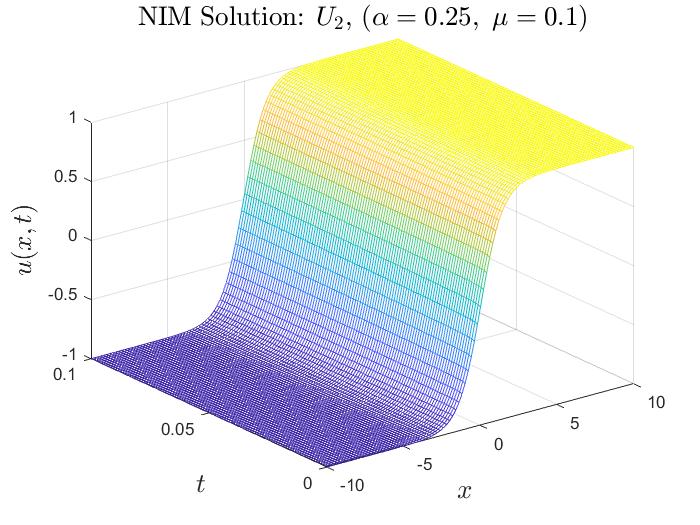}
        %\caption{}
        %\label{Fig1u}
    \end{subfigure}
    ~ %add desired spacing between images, e. g. ~, \quad, \qquad, \hfill etc. 
    \begin{subfigure}[b]{0.40\textwidth}
        \includegraphics[width=\textwidth]{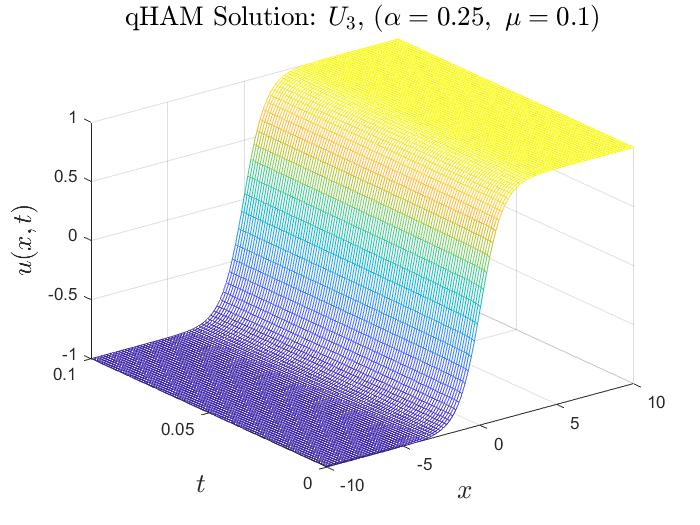}
        %\caption{}
        %\label{fi2}
    \end{subfigure}
\caption{Case I: Comparison between NIM and q-HAM Solution for $n=1$ and $h=-1.$}\label{fig:Fig4}
\end{figure}

\begin{figure}[H]
    \centering
    \begin{subfigure}[b]{0.40\textwidth}
        \includegraphics[width=\textwidth]{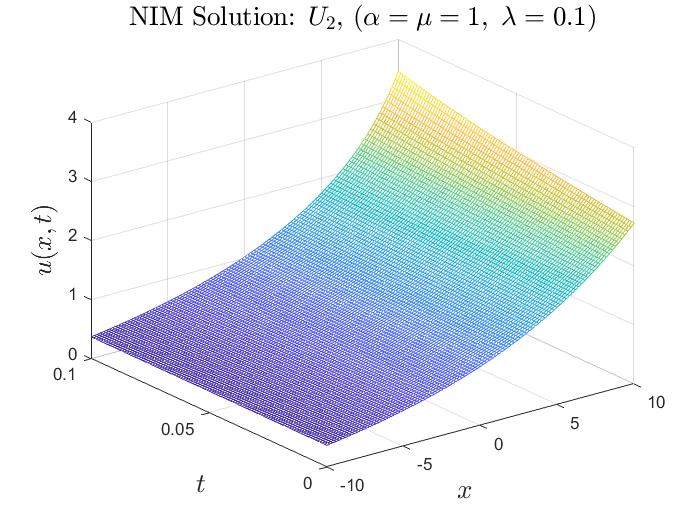}
        %\caption{}
        %\label{Fig1u}
    \end{subfigure}
    ~ %add desired spacing between images, e. g. ~, \quad, \qquad, \hfill etc. 
    \begin{subfigure}[b]{0.40\textwidth}
        \includegraphics[width=\textwidth]{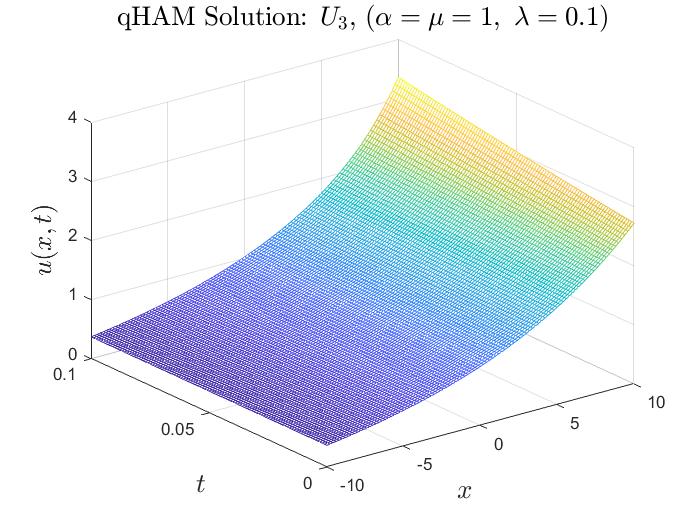}
        %\caption{}
        %\label{fi2}
    \end{subfigure}
\caption{Case II: Comparison between NIM and q-HAM Solution for $n=1$ and $h=-1.$}\label{fig:Fig5}
\end{figure}

\begin{figure}[H]
    \centering
    \begin{subfigure}[b]{0.40\textwidth}
        \includegraphics[width=\textwidth]{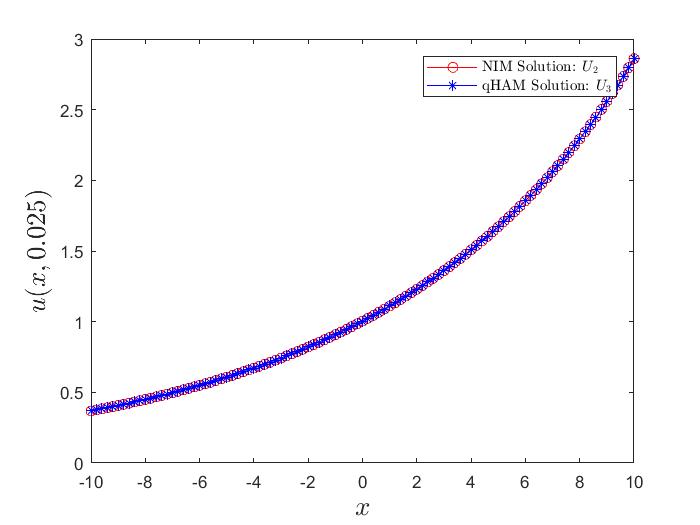}
        %\caption{}
        %\label{Fig1u}
    \end{subfigure}
\caption{Case II: Comparison between NIM and q-HAM Solution in 2D for $n=\alpha=\mu=1,$\ $h=-1,$ and $\lambda=0.1$}\label{fig:Fig6}
\end{figure}

\begin{figure}[H]
    \centering
    \begin{subfigure}[b]{0.40\textwidth}
\includegraphics[width=\textwidth]{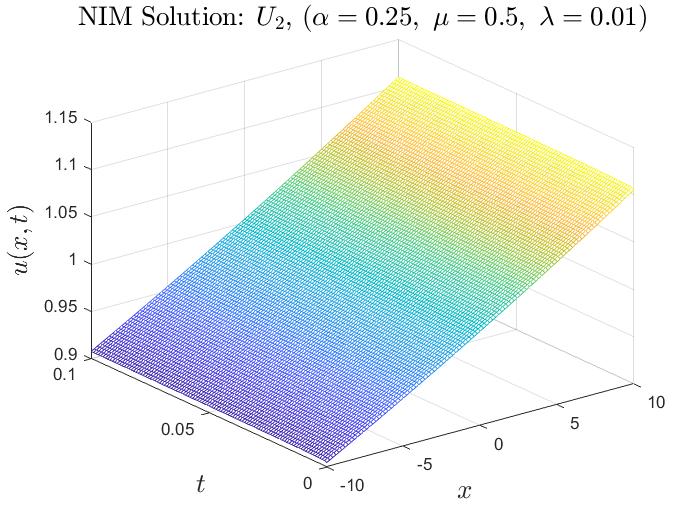}
       % \caption{}
        %\label{Fig1u}
    \end{subfigure}
    ~ %add desired spacing between images, e. g. ~, \quad, \qquad, \hfill etc. 
    \begin{subfigure}[b]{0.40\textwidth}
\includegraphics[width=\textwidth]{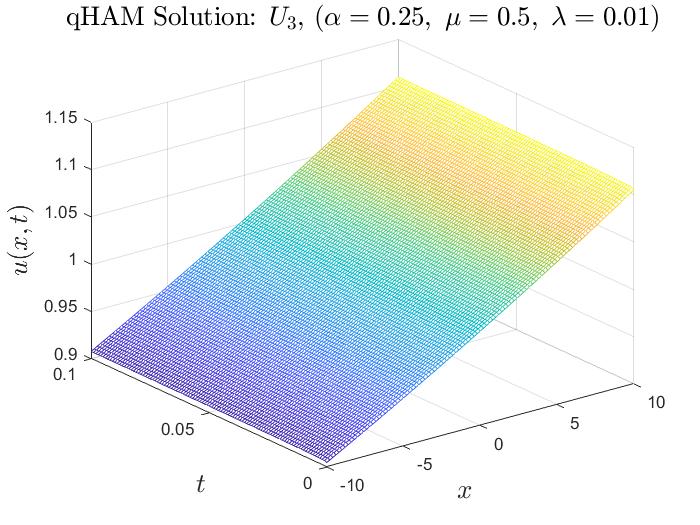}
        %\caption{}
        %\label{fi2}
    \end{subfigure}
\caption{Case II: Comparison between NIM and q-HAM Solution for $n=1$ and $h=-1.$}\label{fig:Fig7}
\end{figure}

Here, we check how accurate these two methods are for solving time-fractional Cahn-Hillard \autoref{eq:1} with different initial conditions as shown in section \ref{sec3} Case I and II. In \Crefrange{fig:Fig1}{fig:Fig7}, one can acknowledge how closely the approximation series solution obtained by these two methods and the exact solution. In \Cref{tab:table1}, error analysis was done for the case when the exact solution is known. 
\\
\\
%%%%%%%%%%%%%%%%%%%%%%%%%%%%%%%%%%%%%%%%%%%%%%%%%%%%%%%%%%%%%%%%
           %%%%%%%%%%%%% Case I Graph %%%%%%%%%%%%%%%%%%
%%%%%%%%%%%%%%%%%%%%%%%%%%%%%%%%%%%%%%%%%%%%%%%%%%%%%%%%%%%%%%%%%%
%%%%%%%%%%%%%%%%%%%%%%%%%%%%%%%%%%%%%%%%%%%%%%%%%%
%%%%%%%%%%%%%%%%%%%%%%%%%%%%%%%%%%%%%%%%%%%%%%%%%%%%%%%%%%%%%%% 
%%%%%%%%%%%%%%%%%%%%%%%%%%%%%%%%%%%%%%%%%%%%%%%%%%%%%%%%%%%%%%% 
%%%%%%%%%%%%%%%%%%%%%%%%%%%%%%%%%%%%%%%%%%%%%%%%%%%%%%%%%%%%%%%%%
             %%%%%%%%%% Case II Graph %%%%%%%%%%%%%%%%%%%%%
%%%%%%%%%%%%%%%%%%%%%%%%%%%%%%%%%%%%%%%%%%%%%%%%%%%%%%%%%%%%%%%%%%%
%%%%%%%%%%%%%%%%%%%%%%%%%%%%%%%%%%%%%%%%%%%%%%%%%%%%%%%%%%%%%%%%%%%
%%%%%%%%%%%%%%%%%%%%%%%%%%%%%%%%%%%%%%%%%%%%%%%%%%%%%
%%%%%%%%%%%%%%%%%%%%%%%%%%%%%%%%%%%%%%%%%%%%%%%%%%%%% 
%%%%%%%%%%%%%%%%%%%%%%%%%%%%%%%%%%%%%%%%%%%%%%%%%%%%%%%%%%%%%%%%%
      %%%%%%%%%%%%% Table 1 Case I %%%%%%%%%%%%%%%%%%
%%%%%%%%%%%%%%%%%%%%%%%%%%%%%%%%%%%%%%%%%%%%%%%%%%%%%%%%%%%%%%%%%
\begin{table}[H]
\caption{Case I: Error Analysis for NIM ($U_2$) and q-HAM ($U_3$) using the Exact Solution ($n=\alpha=\mu=1,$\ $h=-1$).}
\centering % used for centering table
\begin{tabular}{c c c c c c c c c} % centered columns (3 columns)
\toprule[1.5pt]
 & &  \multicolumn{2}{c}{Absolute Error} && & & \multicolumn{2}{c}{Absolute Error} \\
\cline{3-4} \cline{8-9}
$t$ & $x$  & NIM  & q-HAM && $t$ & $x$ & NIM  & q-HAM \\
\toprule[1.5pt] \\
0.01 & 0.0 & $1.151971\times 10^{-7}$ & $2.356975\times 10^{-12}$ && 0.05 & 0.0 & $1.306675\times 10^{-5}$ & $ 7.361971\times 10^{-9}$  \\ [0.5ex]
\cline{2-4} \cline{7-9} \\
    & 1.0 & $ 1.810671\times 10^{-7}$ & $ 2.823765\times 10^{-10}$ &&  & 1.0 & $2.224480\times 10^{-5}$ & $1.736922\times 10^{-7}$\\ [0.5ex]
\cline{2-4} \cline{7-9} \\ 
    & 2.0 & $6.167394\times 10^{-8}$ & $5.749512\times 10^{-11}$ &&  & 2.0 &$ 7.794449\times 10^{-6}$ & $ 3.622408\times 10^{-8}$\\ [0.5ex]
 \cline{2-4} \cline{7-9} \\
    & 3.0 & $1.165205\times 10^{-9}$ & $3.757261\times 10^{-11}$ &&  & 3.0 &$1.257660\times 10^{-7}$ & $ 2.328496\times 10^{-8}$\\ [0.5ex]
 %\toprule[1.5pt]
 \toprule[1.5pt] \\
0.08& 0.0 & $4.940148\times 10^{-5}$ & $7.713501\times 10^{-8}$ && 0.1 & 0.0 & $9.109940\times 10^{-5}$ & $2.352262\times 10^{-7}$  \\ [0.5ex]
\cline{2-4} \cline{7-9} \\
    & 1.0 & $ 8.990891\times 10^{-5}$ & $ 1.124520\times 10^{-6}$ &&  & 1.0 & $1.740220\times 10^{-4}$ & $2.722916\times 10^{-6}$\\ [0.5ex]
\cline{2-4} \cline{7-9} \\ 
    & 2.0 & $3.218897\times 10^{-5}$ & $2.387229\times 10^{-7}$ &&  & 2.0 &$ 6.321236\times 10^{-5}$ & $ 5.848640\times 10^{-7}$\\ [0.5ex]
 \cline{2-4} \cline{7-9} \\
    & 3.0 & $4.548965\times 10^{-7}$ & $1.516340\times 10^{-7}$ &&  & 3.0 &$8.108096\times 10^{-7}$ & $ 3.686350\times 10^{-7}$\\ [0.5ex]
 \toprule[1.5pt]
\end{tabular}\label{tab:table1} 
% is used to refer this table in the text
\end{table}
%%%%%%%%%%%%%%%%%%%%%%%%%%%%%%%%%%%%%%%%%%%%%%%%%%%%%%%%%%%%%%%%%
\section{Solutions of Sixth-Order TFCH Equation}
In this section, we presents the application of the above mentioned methods to obtain approximate solution of the sixth-order TFCH \autoref{eq:1N} subject to different initial conditions.
%%%%%%%%%%%%%%%%%%%%%%%%%%%%%%%%%%%%%%%%%%%%%%%%%%%%%%%%%%%%%%%%%
     %%%%%%%%%%%%%%%%%%%%%% CASE I %%%%%%%%%%%%%%%%%%%%%%%%%
%%%%%%%%%%%%%%%%%%%%%%%%%%%%%%%%%%%%%%%%%%%%%%%%%%%%%%%%%%%%%%%%%
\subsection{}
\textbf{Case I:} Consider the TFCH equation \autoref{eq:1N},
\begin{eqnarray}\label{eq:cahn1}
\capu_t^{\alpha}u=\mu uu_x-18uu^2_{xx}-36u^2_xu_{xx}-24uu_{x}u_{xxx}-3u^2u_{xxxx}+u_{xxxx}+u_{xxxxxx}, \quad 0<\alpha \leq 1,
\end{eqnarray}
with the initial condition
\begin{eqnarray}\label{eq:cahn1IC}
u(x,0)=f(x)=\tanh{\left(\frac{x}{\sqrt{2}}\right)}.
\end{eqnarray}
%%%%%%%%%%%%%%%%%%%%%%%%%%%%%%%%%%%%%%%%%%%%%%%%%%%%%%%%%%%%%%%%%
   %%%%%%%%%%% NIM %%%%%%%%%%
%%%%%%%%%%%%%%%%%%%%%%%%%%%%%%%%%%%%%%%%%%%%%%%%%%%%%%%%%%%%%%%%%
\textbf{NIM solution:}\\
Applying $J^\alpha$ to both sides of \autoref{eq:cahn1}, then the IVP \ref{eq:cahn1}-\ref{eq:cahn1IC} is equivalent to the integral equation:
\begin{eqnarray*}
u(x, t)&=& f(x, t)+\mathscr{L}(u(x, t))+\mathcal{N}(u(x, t)),
\end{eqnarray*}
where, 
\begin{eqnarray*}
u_0&=&f(x)=\tanh{\left(\frac{x}{\sqrt{2}}\right)},\\
\mathscr{L}(u)&=&J^\alpha\left(u_{xxxxxx}+u_{xxxx}\right),\\
\mathcal{N}(u)&=&J^\alpha\left(\mu uu_x-18uu^2_{xx}-36u^2_xu_{xx}-24uu_{x}u_{xxx}-3u^2u_{xxxx}\right).
\end{eqnarray*}
We now obtain components of the series solution using NIM recurrent relation in \autoref{eq:func6} successively as follows:
\begin{eqnarray*}
u_1(x,t)&=& \mathscr{L}(u_0)+\mathcal{N}(u_0)\\
&=& J^\alpha\left(u_{0xxxxxx}+u_{0xxxx}+\mu u_0u_{0x}-18u_0u^2_{0xx}\right)\\
&=& J^\alpha\left(-36u^2_{0x}u_{0xx}-24u_0u_{0x}u_{0xxx}-3u^2_{0}u_{0xxxx}\right)\\
&=& \frac{\mu \tanh{\left(\frac{x}{\sqrt{2}}\right)}\sech^2{\left(\frac{x}{\sqrt{2}}\right)}}{\sqrt{2}\Gamma(\alpha+1)}t^{\alpha},
\end{eqnarray*}

\begin{eqnarray*}
u_2(x,t)&=& \mathscr{L}(u_1)+\big\{\mathcal{N}(u_0+u_1)-\mathcal{N}(u_0)\big\}\\
&=& J^\alpha\left(u_{1xxxxxx}+u_{1xxxx}-36u^2_{0x}u_{1xx}-72u_{0x}u_{1x}u_{0xx}-24u_{0}u_{0x}u_{1xxx}-24u_{0}u_{1x}u_{0xxx}\right)\\
&&-J^\alpha\left(24u_{1}u_{0x}u_{0xxx}+36u_{0}u_{0xx}u_{1xx}+18u_{1}u^2_{0xx}+6u_{0}u_{1}u_{0xxxx}+3u^2_{0}u_{1xxxx}-\mu u_0u_{1}-\mu u_1u_{0}\right)\\
&&-J^\alpha\left(72u_{0x}u_{1x}u_{1xx}+36u^2_{1x}u_{0xx}+24u_{0}u_{1x}u_{1xxx}+24u_{1}u_{0x}u_{1xxx}+24u_{1}u_{1x}u_{0xxx}+18u_{0}u^2_{1xx}-\mu u_{1}u_{1x}\right)\\
&&-J^\alpha\left(36u_{1x}u_{0xx}u_{1xx}+3u^2_{1}u_{0xxxx}+6u_{0}u_{1}u_{1xxxx}+36u^2_{1x}u^2_{1xx}+24u_{1}u_{1x}u_{1xxx}+18u_{1}u^2_{1xx}+3u^2_{1}u_{1xxxx}\right)\\
&=& -\frac{\mu \tanh{\left(\frac{x}{\sqrt{2}}\right)}\sech^8{\left(\frac{x}{\sqrt{2}}\right)}\left(\mu \cosh^6{\left(\frac{x}{\sqrt{2}}\right)}+\Big(96\sqrt{2} -2\mu\Big)\cosh^4{\left(\frac{x}{\sqrt{2}}\right)}-585\sqrt{2}\cosh^2{\left(\frac{x}{\sqrt{2}}\right)}+630\sqrt{2})\right)}{\Gamma(2\alpha+1)}t^{2\alpha}\\
&&+\frac{\mu^2\Gamma(2\alpha+1)\tanh{\left(\frac{x}{\sqrt{2}}\right)}\sech^{4}{\left(\frac{x}{\sqrt{2}}\right)}\left(3\Big(\mu\sqrt{2}+1428\Big)\sech^{2}{\left(\frac{x}{\sqrt{2}}\right)}-2\Big(\sqrt{2}\mu+192\Big)\right)}{64\Gamma(\alpha+1)^2\Gamma(3\alpha+1)}t^{3\alpha}\\
&&+\frac{420\mu^2\Gamma(2\alpha+1)\tanh{\left(\frac{x}{\sqrt{2}}\right)}\sech^{10}{\left(\frac{x}{\sqrt{2}}\right)}\Big(5-13\cosh(\sqrt{2}x)\Big)}{64\Gamma(\alpha+1)^2\Gamma(3\alpha+1)}t^{3\alpha}\\
&&-\frac{3\mu^3\Gamma(3\alpha+1)\tanh{\left(\frac{x}{\sqrt{2}}\right)}\sech^{12}{\left(\frac{x}{\sqrt{2}}\right)}\left(3773\cosh(\sqrt{2}x)-646\cosh(2\sqrt{2}x)+27\cosh(3\sqrt{2}x)-3474\right)}{16\sqrt{2}\Gamma(\alpha+1)^3\Gamma(4\alpha+1)}t^{4\alpha}. 
\end{eqnarray*}
Using the same procedure, expression for $u_m(x,t),$ $m=3,4,5,...$ can be obtained. The expression of the series solutions given by this new iteration (NIM) can be written in the form

\begin{eqnarray}\notag
u(x,t)&\cong&\sum_{i=0}^{2} u_i(x,t)=u_0(x,t)+u_1(x,t)+u_2(x,t)\\ \notag
&=& \tanh{\left(\frac{x}{\sqrt{2}}\right)}+\frac{\mu \tanh{\left(\frac{x}{\sqrt{2}}\right)}\sech^2{\left(\frac{x}{\sqrt{2}}\right)}}{\sqrt{2}\Gamma(\alpha+1)}t^\alpha\\ \notag
&&-\frac{\mu \tanh{\left(\frac{x}{\sqrt{2}}\right)}\sech^8{\left(\frac{x}{\sqrt{2}}\right)}\left(\mu \cosh^6{\left(\frac{x}{\sqrt{2}}\right)}+\Big(96\sqrt{2} -2\mu\Big)\cosh^4{\left(\frac{x}{\sqrt{2}}\right)}-585\sqrt{2}\cosh^2{\left(\frac{x}{\sqrt{2}}\right)}+630\sqrt{2})\right)}{\Gamma(2\alpha+1)}t^{2\alpha}\\ \notag
&&+\frac{\mu^2\Gamma(2\alpha+1)\tanh{\left(\frac{x}{\sqrt{2}}\right)}\sech^{4}{\left(\frac{x}{\sqrt{2}}\right)}\left(3\Big(\mu\sqrt{2}+1428\Big)\sech^{2}{\left(\frac{x}{\sqrt{2}}\right)}-2\Big(\sqrt{2}\mu+192\Big)\right)}{64\Gamma(\alpha+1)^2\Gamma(3\alpha+1)}t^{3\alpha}\\
&&+\frac{420\mu^2\Gamma(2\alpha+1)\tanh{\left(\frac{x}{\sqrt{2}}\right)}\sech^{10}{\left(\frac{x}{\sqrt{2}}\right)}\Big(5-13\cosh(\sqrt{2}x)\Big)}{64\Gamma(\alpha+1)^2\Gamma(3\alpha+1)}t^{3\alpha}\\ \label{eq:solutcahn}
&&-\frac{3\mu^3\Gamma(3\alpha+1)\tanh{\left(\frac{x}{\sqrt{2}}\right)}\sech^{12}{\left(\frac{x}{\sqrt{2}}\right)}\left(3773\cosh(\sqrt{2}x)-646\cosh(2\sqrt{2}x)+27\cosh(3\sqrt{2}x)-3474\right)}{16\sqrt{2}\Gamma(\alpha+1)^3\Gamma(4\alpha+1)}t^{4\alpha}. \notag
\end{eqnarray}
Thus, \autoref{eq:solutcahn} gives an approximate solution to the IVP \ref{eq:cahn1}-\ref{eq:cahn1IC}.\\
%%%%%%%%%%%%%%%%%%%%%%%%%%%%%%%%%%%%%%%%%%%%%%%%%%%%%%%%%%%%%%%%%
        %%%%%%%%%%%%%%%% q-HAM %%%%%%%%%%%%%%%%%%%
%%%%%%%%%%%%%%%%%%%%%%%%%%%%%%%%%%%%%%%%%%%%%%%%%%%%%%%%%%%%%%%%%
\newline
\textbf{q-HAM solution:}\\
To apply the q-HAM, we rewrite TFCH \autoref{eq:1} as, 
\begin{eqnarray}\label{eq:cahn3}
\capu_t^{\alpha}u-\mu uu_x+18uu^2_{xx}+36u^2_xu_{xx}+24uu_{x}u_{xxx}+3u^2u_{xxxx}-u_{xxxx}-u_{xxxxxx}=0,\quad 0<\alpha \leq 1. \end{eqnarray}
 Applying q-HAM to \autoref{eq:cahn3}, we obtain from \autoref{eq:8a},
\begin{eqnarray}\notag
\mathcal{R}_m\left(\vec{u}_{m-1}\right)&=&\capu_t^\alpha u_{(m-1)}-\mu\sum_{k=0}^{m-1} u_{k}u_{(m-1-k)x}+18\sum_{k=0}^{m-1}\sum_{j=0}^{k}u_{j} u_{(k-j)xx}u_{(m-1-k)xx}\\\notag
&& +36\sum_{k=0}^{m-1}\sum_{j=0}^{k}u_{jx} u_{(k-j)x}u_{(m-1-k)xx}+24\sum_{k=0}^{m-1}\sum_{j=0}^{k}u_j u_{(k-j)x}u_{(m-1-k)xxx}\\ 
&& +3\sum_{k=0}^{m-1}\sum_{j=0}^{k}u_{j} u_{(k-j)}u_{(m-1-k)xxxx} -u_{(m-1)xxxx}-u_{(m-1)xxxxxx}.\label{eq:cahn6s}
\end{eqnarray}
We can now obtain the components of the solution using q-HAM recurrent relation in \autoref{eq:6ss}, using Eqs.~\ref{eq:9a} and \ref{eq:cahn6s} successively as follows:
\begin{eqnarray*}
u_1(x,t)&=&\Psi^*_1 u_0(x,t)+h \mathcal{R}_1\left(\vec{u}_0(x,t)\right)\\
&=& h J^\alpha\left(\mathcal{D}_t^\alpha u_0 -\mu u_0u_{0x}+18u_0u_{0xx}u_{0xx}+36u_{0x}u_{0x}u_{0xx}\right)\\
&&+h J^\alpha\left(24u_0u_{0x}u_{0xxx}+3u_0u_{0}u_{0xxxx}-u_{0xxxx}-u_{0xxxxxx}\right)\\
&=& -\frac{h\mu \tanh{\left(\frac{x}{\sqrt{2}}\right)}\sech^2{\left(\frac{x}{\sqrt{2}}\right)}}{\sqrt{2}\Gamma(\alpha+1)}t^{\alpha},
\end{eqnarray*}

\begin{eqnarray*}
u_2(x,t)&=&\Psi^*_2 u_1(x,t)+h J^\alpha\left(\mathcal{R}_2\left(\vec{u}_1(x,t)\right)\right)\\
&=&hu_1+ h J^\alpha\left(\mathcal{D}_t^\alpha u_1-\mu u_{0}u_{1x}-\mu u_{1}u_{0x}+18u_0u_{0xx}u_{1xx}+18u_0u_{1xx}u_{0xx}+18u_1u_{0xx}u_{0xx}\right)\\
&&+ \ h J^\alpha\left(36u_{0x}u_{0x}u_{1xx}+36u_{0x}u_{1x}u_{0xx}+36u_{1x}u_{0x}u_{0xx}+24u_0u_{0x}u_{1xxx}+24u_0u_{1x}u_{0xxx}\right)\\
&&+ \ h J^\alpha\left(24u_{1}u_{0x}u_{0xxx}+3u_{0}u_{0}u_{1xxxx}+3u_{0}u_{1}u_{0xxxx}+3u_{1}u_{0}u_{0xxxx}-u_{1xxxx}-u_{1xxxxxx}\right)\\
&=& (n+h)\,u_1-\frac{h^2\mu \tanh{\left(\frac{x}{\sqrt{2}}\right)}\sech^8{\left(\frac{x}{\sqrt{2}}\right)}\left(\mu \cosh^6{\left(\frac{x}{\sqrt{2}}\right)}+\Big(96\sqrt{2} -2\mu\Big)\cosh^4{\left(\frac{x}{\sqrt{2}}\right)}\right)}{\Gamma(2\alpha+1)}t^{2\alpha}\\
&&-\frac{h^2\mu \tanh{\left(\frac{x}{\sqrt{2}}\right)}\sech^8{\left(\frac{x}{\sqrt{2}}\right)}\left(
-585\sqrt{2}\cosh^2{\left(\frac{x}{\sqrt{2}}\right)}+630\sqrt{2})\right)}{\Gamma(2\alpha+1)}t^{2\alpha},
\end{eqnarray*}

\begin{eqnarray*}
u_3(x,t)&=&\Psi^*_3 u_2(x,t)+h J^\alpha\left(\mathcal{R}_3\left(\vec{u}_2(x,t)\right)\right)\\
&=&hu_2+ h J^\alpha\left(\mathcal{D}_t^\alpha u_2-\mu u_{0}u_{2x}-\mu u_{1}u_{1x}-\mu u_{2}u_{0x}+18u_0u_{0xx}u_{2xx}+18u_0u_{1xx}u_{1xx}+18u_1u_{0xx}u_{1xx}\right)\\
&&+ \ h J^\alpha\left(18u_0u_{2xx}u_{0xx}+18u_1u_{1xx}u_{0xx}+18u_2u_{0xx}u_{0xx}+36u_{0x}u_{0x}u_{2xx}+36u_{0x}u_{1x}u_{1xx}+36u_{1x}u_{0x}u_{1xx}\right)\\
&&+ \ h J^\alpha\left(36u_{0x}u_{2x}u_{0xx}+36u_{1x}u_{0x}u_{1xx}+36u_{2x}u_{0x}u_{0xx}+24u_0u_{0x}u_{2xxx}+24u_0u_{1x}u_{1xxx}+24u_1u_{0x}u_{1xxx}\right)\\
&&+ \ h J^\alpha\left(24u_0u_{2x}u_{0xxx}+24u_0u_{1x}u_{1xxx}+24u_2u_{0x}u_{0xxx}+3u_{0}u_{0}u_{2xxxx}+3u_{0}u_{1}u_{1xxxx}+3u_{1}u_{0}u_{1xxxx}\right)\\
&&+ \ h J^\alpha\left(3u_{0}u_{2}u_{0xxxx}+3u_{1}u_{1}u_{0xxxx}+3u_{2}u_{0}u_{0xxxx}-u_{2xxxx}-u_{2xxxxxx}\right)\\
&=& (n+h)\,u_2-\frac{\mu h^2(n+h) \tanh{\left(\frac{x}{\sqrt{2}}\right)}\sech^8{\left(\frac{x}{\sqrt{2}}\right)}\left(\mu \cosh^6{\left(\frac{x}{\sqrt{2}}\right)}+\Big(96\sqrt{2} -2\mu\Big)\cosh^4{\left(\frac{x}{\sqrt{2}}\right)}\right)}{\Gamma(2\alpha+1)}t^{2\alpha}\\
&&-\frac{\mu h^2(n+h) \tanh{\left(\frac{x}{\sqrt{2}}\right)}\sech^8{\left(\frac{x}{\sqrt{2}}\right)}\left(-585\sqrt{2}\cosh^2{\left(\frac{x}{\sqrt{2}}\right)}+630\sqrt{2})\right)}{\Gamma(2\alpha+1)}t^{2\alpha}\\
&&+\frac{144\sqrt{2}h^3\mu \tanh\left(\frac{x}{\sqrt{2}}\right) \sech^{10}{\left(\frac{x}{\sqrt{2}}\right)}\Big(4484\cosh(\sqrt{2}x)-471\cosh(2\sqrt{2}x)+8\cosh(3\sqrt{2}x)-5117\Big)}{\Gamma(3\alpha+1)}t^{3\alpha}\\
&&-\frac{3h^3\mu^2 \tanh\left(\frac{x}{\sqrt{2}}\right) \sech^{10}{\left(\frac{x}{\sqrt{2}}\right)}\Big(17972\cosh(\sqrt{2}x)-2031\cosh(2\sqrt{2}x)+56\cosh(3\sqrt{2}x)-20261\Big)}{4\Gamma(3\alpha+1)}t^{3\alpha}\\
&&+\frac{h^3\mu^3 \tanh\left(\frac{x}{\sqrt{2}}\right) \sech^{6}{\left(\frac{x}{\sqrt{2}}\right)}\Big(-20\cosh(\sqrt{2}x)+\cosh(2\sqrt{2}x)+27\Big)}{4\sqrt{2}\Gamma(3\alpha+1)}t^{3\alpha}\\%%
&&+\frac{3h^3\mu^2\Gamma(2\alpha+1) \tanh\left(\frac{x}{\sqrt{2}}\right) \sech^{10}{\left(\frac{x}{\sqrt{2}}\right)}\Big(2323\cosh(\sqrt{2}x)-309\cosh(2\sqrt{2}x)+8\cosh(3\sqrt{2}x)-2391\Big)}{8\Gamma(\alpha+1)^2\Gamma(3\alpha+1)}t^{3\alpha}\\
&&+\frac{h^3\mu^3\Gamma(2\alpha+1) \tanh\left(\frac{x}{\sqrt{2}}\right) \sech^{6}{\left(\frac{x}{\sqrt{2}}\right)}\Big(\cosh(\sqrt{2}x)-2\Big)}{2\sqrt{2}\Gamma(\alpha+1)^2\Gamma(3\alpha+1)}t^{3\alpha}\\
&&+\frac{9h^3\mu \tanh\left(\frac{x}{\sqrt{2}}\right) \sech^{14}{\left(\frac{x}{\sqrt{2}}\right)}\Big(153385672396573417792\cosh\left(\sqrt{2}x\right)- 216241588341872290245\Big)}{17592186044416\Gamma(3\alpha+1)}t^{3\alpha}\\
&&-\frac{9h^3\mu \tanh\left(\frac{x}{\sqrt{2}}\right) \sech^{14}{\left(\frac{x}{\sqrt{2}}\right)}\Big(69679073538343302995\cosh\left(2\sqrt{2}x\right)\Big)}{17592186044416\Gamma(3\alpha+1)}t^{3\alpha}\\%
&&-\frac{3h^3\mu \tanh\left(\frac{x}{\sqrt{2}}\right) \sech^{8}{\left(\frac{x}{\sqrt{2}}\right)}\Big(49120689046652631\cosh\left(\sqrt{2}x\right)-1487475224720429629\Big)}{68719476736\Gamma(3\alpha+1)}t^{3\alpha}\\
&&+\frac{h^3\mu^2\tanh\left(\frac{x}{\sqrt{2}}\right)\sech^{8}{\left(\frac{x}{\sqrt{2}}\right)}\Big(5457854338516959\cosh\left(2\sqrt{2}x\right)-24770261997884659\Big)
}{2199023255552\sqrt{2}\Gamma(3\alpha+1)}t^{3\alpha}\\
&&-\frac{39184595250890985h^3\mu^2\tanh\left(\frac{x}{\sqrt{2}}\right)\sech^{10}{\left(\frac{x}{\sqrt{2}}\right)}
}{4398046511104\sqrt{2}\Gamma(3\alpha+1)}t^{3\alpha}.
\end{eqnarray*}

Using the same procedure, expression for $u_m(x,t),$ $m=4,5,6,...$ can be obtained. The expression of the series solutions given by q-HAM can be written in the form

\begin{eqnarray}\notag
u(x,t;n;h)&\cong&u_0(x,t)+\sum_{i=1}^3 u_i(x,t;n;h)\left(\frac{1}{n}\right)^i\\ \notag
&=& \tanh{\left(\frac{x}{\sqrt{2}}\right)}-\frac{\mu h(3n^2+3nh+h^2)\tanh{\left(\frac{x}{\sqrt{2}}\right)}\sech^2{\left(\frac{x}{\sqrt{2}}\right)}}{n^3\Gamma(\alpha+1)}t^{\alpha}\\ \notag
&&-\frac{\mu h^2(3n+2h)\tanh{\left(\frac{x}{\sqrt{2}}\right)}\sech^8{\left(\frac{x}{\sqrt{2}}\right)}\left(\mu \cosh^6{\left(\frac{x}{\sqrt{2}}\right)}+\Big(96\sqrt{2} -2\mu\Big)\cosh^4{\left(\frac{x}{\sqrt{2}}\right)}\right)}{n^3\Gamma(2\alpha+1)}t^{2\alpha}\\ \notag
&&-\frac{\mu h^2(3n+2h) \tanh{\left(\frac{x}{\sqrt{2}}\right)}\sech^8{\left(\frac{x}{\sqrt{2}}\right)}\left(
-585\sqrt{2}\cosh^2{\left(\frac{x}{\sqrt{2}}\right)}+630\sqrt{2})\right)}{n^3\Gamma(2\alpha+1)}t^{2\alpha}\\ \notag%%%
&&+\frac{144\sqrt{2}h^3\mu \tanh\left(\frac{x}{\sqrt{2}}\right) \sech^{10}{\left(\frac{x}{\sqrt{2}}\right)}\Big(4484\cosh(\sqrt{2}x)-471\cosh(2\sqrt{2}x)+8\cosh(3\sqrt{2}x)-5117\Big)}{n^3\Gamma(3\alpha+1)}t^{3\alpha}\\ \notag
&&-\frac{3h^3\mu^2 \tanh\left(\frac{x}{\sqrt{2}}\right) \sech^{10}{\left(\frac{x}{\sqrt{2}}\right)}\Big(17972\cosh(\sqrt{2}x)-2031\cosh(2\sqrt{2}x)+56\cosh(3\sqrt{2}x)-20261\Big)}{4n^3\Gamma(3\alpha+1)}t^{3\alpha}\\ \notag
&&+\frac{h^3\mu^3 \tanh\left(\frac{x}{\sqrt{2}}\right) \sech^{6}{\left(\frac{x}{\sqrt{2}}\right)}\Big(-20\cosh(\sqrt{2}x)+\cosh(2\sqrt{2}x)+27\Big)}{4\sqrt{2}n^3\Gamma(3\alpha+1)}t^{3\alpha}\\ \notag
&&+\frac{3h^3\mu^2\Gamma(2\alpha+1) \tanh\left(\frac{x}{\sqrt{2}}\right) \sech^{10}{\left(\frac{x}{\sqrt{2}}\right)}\Big(2323\cosh(\sqrt{2}x)-309\cosh(2\sqrt{2}x)+8\cosh(3\sqrt{2}x)-2391\Big)}{8n^3\Gamma(\alpha+1)^2\Gamma(3\alpha+1)}t^{3\alpha}\\ \notag
&&+\frac{h^3\mu^3\Gamma(2\alpha+1) \tanh\left(\frac{x}{\sqrt{2}}\right) \sech^{6}{\left(\frac{x}{\sqrt{2}}\right)}\Big(\cosh(\sqrt{2}x)-2\Big)}{2\sqrt{2}n^3\Gamma(\alpha+1)^2\Gamma(3\alpha+1)}t^{3\alpha}\\ \notag
&&+\frac{9h^3\mu \tanh\left(\frac{x}{\sqrt{2}}\right) \sech^{14}{\left(\frac{x}{\sqrt{2}}\right)}\Big(153385672396573417792\cosh\left(\sqrt{2}x\right)- 216241588341872290245\Big)}{17592186044416n^3\Gamma(3\alpha+1)}t^{3\alpha}\\ \notag
&&-\frac{9h^3\mu \tanh\left(\frac{x}{\sqrt{2}}\right) \sech^{14}{\left(\frac{x}{\sqrt{2}}\right)}\Big(69679073538343302995\cosh\left(2\sqrt{2}x\right)\Big)}{17592186044416n^3\Gamma(3\alpha+1)}t^{3\alpha}\\ \notag
&&-\frac{3h^3\mu \tanh\left(\frac{x}{\sqrt{2}}\right) \sech^{8}{\left(\frac{x}{\sqrt{2}}\right)}\Big(49120689046652631\cosh\left(\sqrt{2}x\right)-1487475224720429629\Big)}{68719476736n^3\Gamma(3\alpha+1)}t^{3\alpha}\\ \notag
&&+\frac{h^3\mu^2\tanh\left(\frac{x}{\sqrt{2}}\right)\sech^{8}{\left(\frac{x}{\sqrt{2}}\right)}\Big(5457854338516959\cosh\left(2\sqrt{2}x\right)-24770261997884659\Big)
}{2199023255552\sqrt{2}n^3\Gamma(3\alpha+1)}t^{3\alpha}\\ 
&&-\frac{39184595250890985h^3\mu^2\tanh\left(\frac{x}{\sqrt{2}}\right)\sech^{10}{\left(\frac{x}{\sqrt{2}}\right)}}{4398046511104\sqrt{2}n^3\Gamma(3\alpha+1)}t^{3\alpha}.\label{eq:solutcahn3}
\end{eqnarray}
Thus, \autoref{eq:solutcahn3} gives an approximate solutions to the IVP \ref{eq:cahn1}-\ref{eq:cahn1IC} in terms of convergence parameter $h$ and $n$.\\
\\
%%%%%%%%%%%%%%%%%%%%%%%%%%%%%%%%%%%%%%%%%%%%%%%%%%%%%%%%%%%%%%%%%
     %%%%%%%%%%%%%%%%%%%%%% CASE II %%%%%%%%%%%%%%%%%%%%%%%%%
%%%%%%%%%%%%%%%%%%%%%%%%%%%%%%%%%%%%%%%%%%%%%%%%%%%%%%%%%%%%%%%%%
\textbf{Case II:} Consider the TFCH \autoref{eq:1N},
\begin{eqnarray}\label{eq:cahn2}
\capu_t^{\alpha}u=\mu uu_x-18uu^2_{xx}-36u^2_xu{xx}-24uu_{x}u_{xxx}-3u^2u_{xxxx}+u_{xxxx}+u_{xxxxxx}, \quad 0<\alpha \leq 1,
\end{eqnarray}
with the initial condition
\begin{eqnarray}\label{eq:cahn2IC}
u(x,0)=f(x)=e^{ \lambda x}.
\end{eqnarray}
%%%%%%%%%%%%%%%%%%%%%%%%%%%%%%%%%%%%%%%%%%%%%%%%%%%%%%%%%%%%%%%%%
%%%%%%%%%%%%%%%%%%%%%%%%%%%%%%%%%%%%%%%%%%%%%%%%%%%%%%%%%%%%%%%%%
          %%%%%%%%%%%%%%%%%% NIM %%%%%%%%%%%%%%%%%%
%%%%%%%%%%%%%%%%%%%%%%%%%%%%%%%%%%%%%%%%%%%%%%%%%%%%%%%%%%%%%%%%%
\newline
\textbf{NIM solution:}\\
Applying $J^\alpha$ to both sides of \autoref{eq:cahn2}, then the IVP \ref{eq:cahn2}-\ref{eq:cahn2IC} is equivalent to the integral equation:
\begin{eqnarray*}
u(x, t)&=& f(x,t)+\mathscr{L}(u(x, t))+\mathcal{N}(u(x,t)),
\end{eqnarray*}
where, 
\begin{eqnarray*}
u_0&=&f(x)=e^{\lambda x},\\
\mathscr{L}(u)&=&J^\alpha\left(u_{xxxxxx}+u_{xxxx}\right),\\
\mathcal{N}(u)&=&J^\alpha\left(\mu uu_x-18uu^2_{xx}-36(u_x)^2u_{xx}-24uu_{x}u_{xxx}-3u^2u_{xxxx}\right).
\end{eqnarray*}
We now obtain components of the series solution using NIM recurrent relation in \autoref{eq:func6} successively as follows:
\begin{eqnarray*}
u_1(x,t)&=& \mathscr{L}(u_0)+\mathcal{N}(u_0)\\
&=& J^\alpha\left(u_{0xxxxxx}+u_{0xxxx}+\mu u_0u_{0x}-18u_0u^2_{0xx}-36(u_{0x})^2u_{0xx}-24u_0u_{0x}u_{0xxx}-3u^2_{0}u_{0xxxx}\right)\\
&=& \frac{\lambda e^{\lambda x}(\lambda^5+\lambda^3(1-81e^{2\lambda x})+\mu e^{\lambda x}) }{\Gamma(\alpha+1)}t^{\alpha},
\end{eqnarray*}

\begin{eqnarray*}
u_2(x,t)&=& J^\alpha\left(\mathscr{L}(u_1)+\big\{\mathcal{N}(u_0+u_1)-\mathcal{N}(u_0)\big\}\right)\\
&=& J^\alpha\left(u_{1xxxxxx}+u_{1xxxx}-36u^2_{0x}u_{1xx}-72u_{0x}u_{1x}u_{0xx}-24u_{0}u_{0x}u_{1xxx}-24u_{0}u_{1x}u_{0xxx}\right)\\
&&-J^\alpha\left(24u_{1}u_{0x}u_{0xxx}+36u_{0}u_{0xx}u_{1xx}+18u_{1}u^2_{0xx}+6u_{0}u_{1}u_{0xxxx}+3u^2_{0}u_{1xxxx}-\mu u_0u_{1}-\mu u_1u_{0}\right)\\
&&-J^\alpha\left(72u_{0x}u_{1x}u_{1xx}+36u^2_{1x}u_{0xx}+24u_{0}u_{1x}u_{1xxx}+24u_{1}u_{0x}u_{1xxx}+24u_{1}u_{1x}u_{0xxx}+18u_{0}u^2_{1xx}-\mu u_{1}u_{1x}\right)\\
&&-J^\alpha\left(36u_{1x}u_{0xx}u_{1xx}+3u^2_{1}u_{0xxxx}+6u_{0}u_{1}u_{1xxxx}+36u^2_{1x}u^2_{1xx}+24u_{1}u_{1x}u_{1xxx}+18u_{1}u^2_{1xx}+3u^2_{1}u_{1xxxx}\right)\\
&=& \frac{\lambda^2e^{\lambda x}\left(151875\lambda^6 e^{4\lambda x}-1092\mu \lambda^3e^{3\lambda x}-3(324\lambda^6(61\lambda^2+7)-\mu^2)e^{2\lambda x}+(\lambda^2+1)^2\lambda^6+6\mu(11\lambda^2+3)\lambda^3e^{\lambda x}\right)}{\Gamma(2\alpha+1)}t^{2\alpha}\\
&&-\frac{\lambda^3e^{4\lambda x}\Gamma(2\alpha+1)\left(-303750\lambda^{11}e^{\lambda x}+47258883\lambda^9e^{3\lambda x}-649539\mu \lambda^6e^{2\lambda x}+1860\mu \lambda^6-2\mu^3\right)}{\Gamma(\alpha+1)^2\Gamma(3\alpha+1)}t^{3\alpha}\\ 
&&+\frac{\lambda^6e^{2\lambda x}\Gamma(2\alpha+1)\left(-243\lambda^{10}e^{\lambda x}-486\lambda^8e^{\lambda x}+\mu \lambda^7-248\lambda^6e^{\lambda x}+303750\lambda^6e^{3\lambda x}-1860\mu \lambda^5e^{2\lambda x}\right)}{\Gamma(\alpha+1)^2\Gamma(3\alpha+1)}t^{3\alpha}\\
&&+\frac{\lambda^6e^{2\lambda x}\Gamma(2\alpha+1)\left(2\mu \lambda^5+\mu \lambda^3+3\lambda^2\mu^2e^{\lambda x}+3\mu^2e^{2\lambda x}-2280\mu e^{3\lambda x}\right)}{\Gamma(\alpha+1)^2\Gamma(3\alpha+1)}t^{3\alpha}\\ 
&&+\frac{3\lambda^{10}e^{5\lambda x}\Gamma(3\alpha+1)\left(101250\lambda^8- 15752961\lambda^6e^{2\lambda x}+1162261467\lambda^6e^{4\lambda x}+50625\lambda^6\right)}{\Gamma(\alpha+1)^3\Gamma(4\alpha+1)}t^{4\alpha}\\ 
&&+\frac{3\lambda^{10}e^{5\lambda x}\Gamma(3\alpha+1)\left(209952\mu \lambda^5e^{\lambda x}-625\mu^2\lambda^2+194481\mu^2e^{2\lambda x}-625\mu^2\right)}{\Gamma(\alpha+1)^3\Gamma(4\alpha+1)}t^{4\alpha}\\ 
&&-\frac{3\lambda^{7}e^{3\lambda x}\Gamma(3\alpha+1)\left(27\lambda^{15}- 50625\lambda^{13}e^{2\lambda x}+18\lambda^{13}+15752961\lambda^{11}e^{4\lambda x}+81\lambda^{11}+256\mu \lambda^{10}e^{\lambda x}\right)}{\Gamma(\alpha+1)^3\Gamma(4\alpha+1)}t^{4\alpha}\\
&&-\frac{3\lambda^{7}e^{3\lambda x}\Gamma(3\alpha+1)\left(27\lambda^{9}+ 512\mu \lambda^{8}e^{\lambda x}+256\mu \lambda^{6}e^{\lambda x}-209952\mu \lambda^{6}e^{3\lambda x}+26873856\mu\lambda^{6}e^{5\lambda x}+432\mu^3e^{\lambda x}\right)}{\Gamma(\alpha+1)^3\Gamma(4\alpha+1)}t^{4\alpha}.
\end{eqnarray*}
Using the same procedure, expression for $u_m(x,t),$ $m=3,4,5,...$ can be obtained. Thus, the expression of the series solutions given by this new iteration (NIM) can be written in the form

\begin{eqnarray}\notag
u(x,t)&\cong&\sum_{i=0}^{2} u_i(x,t)=u_0(x,t)+u_1(x,t)+u_2(x,t)\\ \notag
&=& e^{\lambda x}+\frac{\lambda e^{\lambda x}(\lambda^5+\lambda^3(1-81e^{2\lambda x})+\mu e^{\lambda x}) }{\Gamma(\alpha+1)}t^\alpha\\ \notag
&&+ \frac{\lambda^2e^{\lambda x}\left(151875\lambda^6 e^{4\lambda x}-1092\mu \lambda^3e^{3\lambda x}-3(324\lambda^6(61\lambda^2+7)-\mu^2)e^{2\lambda x}+(\lambda^2+1)^2\lambda^6+6\mu(11\lambda^2+3)\lambda^3e^{\lambda x}\right)}{\Gamma(2\alpha+1)}t^{2\alpha}\\ \notag
&&-\frac{\lambda^3e^{4\lambda x}\Gamma(2\alpha+1)\left(-303750\lambda^{11}e^{\lambda x}+47258883\lambda^9e^{3\lambda x}-649539\mu \lambda^6e^{2\lambda x}+1860\mu \lambda^6-2\mu^3\right)}{\Gamma(\alpha+1)^2\Gamma(3\alpha+1)}t^{3\alpha}\\ \notag
&&+\frac{\lambda^6e^{2\lambda x}\Gamma(2\alpha+1)\left(-243\lambda^{10}e^{\lambda x}-486\lambda^8e^{\lambda x}+\mu \lambda^7-248\lambda^6e^{\lambda x}+303750\lambda^6e^{3\lambda x}-1860\mu \lambda^5e^{2\lambda x}\right)}{\Gamma(\alpha+1)^2\Gamma(3\alpha+1)}t^{3\alpha}\\ \notag
&&+\frac{\lambda^6e^{2\lambda x}\Gamma(2\alpha+1)\left(2\mu \lambda^5+\mu \lambda^3+3\lambda^2\mu^2e^{\lambda x}+3\mu^2e^{2\lambda x}-2280\mu e^{3\lambda x}\right)}{\Gamma(\alpha+1)^2\Gamma(3\alpha+1)}t^{3\alpha}\\ \notag
&&+\frac{3\lambda^{10}e^{5\lambda x}\Gamma(3\alpha+1)\left(101250\lambda^8- 15752961\lambda^6e^{2\lambda x}+1162261467\lambda^6e^{4\lambda x}+50625\lambda^6\right)}{\Gamma(\alpha+1)^3\Gamma(4\alpha+1)}t^{4\alpha}\\ \notag
&&+\frac{3\lambda^{10}e^{5\lambda x}\Gamma(3\alpha+1)\left(209952\mu \lambda^5e^{\lambda x}-625\mu^2\lambda^2+194481\mu^2e^{2\lambda x}-625\mu^2\right)}{\Gamma(\alpha+1)^3\Gamma(4\alpha+1)}t^{4\alpha}\\ \notag
&&-\frac{3\lambda^{7}e^{3\lambda x}\Gamma(3\alpha+1)\left(27\lambda^{9}+ 512\mu \lambda^{8}e^{\lambda x}+256\mu \lambda^{6}e^{\lambda x}-209952\mu \lambda^{6}e^{3\lambda x}+26873856\mu\lambda^{6}e^{5\lambda x}+432\mu^3e^{\lambda x}\right)}{\Gamma(\alpha+1)^3\Gamma(4\alpha+1)}t^{4\alpha}\\
&&-\frac{3\lambda^{7}e^{3\lambda x}\Gamma(3\alpha+1)\left(27\lambda^{15}- 50625\lambda^{13}e^{2\lambda x}+18\lambda^{13}+15752961\lambda^{11}e^{4\lambda x}+81\lambda^{11}+256\mu \lambda^{10}e^{\lambda x}\right)}{\Gamma(\alpha+1)^3\Gamma(4\alpha+1)}t^{4\alpha}.\label{eq:solutnim1nim6}
\end{eqnarray}
Thus, \autoref{eq:solutnim1nim6} gives an approximate solution to the IVP \ref{eq:cahn2}-\ref{eq:cahn2IC}.\\
%%%%%%%%%%%%%%%%%%%%%%%%%%%%%%%%%%%%%%%%%%%%%%%%%%%%%%%%%%%%%%%%%
          %%%%%%%%%%%%%%%%%% q-HAM %%%%%%%%%%%%%%%%%%
%%%%%%%%%%%%%%%%%%%%%%%%%%%%%%%%%%%%%%%%%%%%%%%%%%%%%%%%%%%%%%%%
\newline
\textbf{q-HAM solution:}\\
Consider \autoref{eq:cahn3}, we obtain the components of the solution using q-HAM recurrent relation in \autoref{eq:6ss}, using Eqs.~\ref{eq:9a} and \ref{eq:cahn6s} successively as follows:
\begin{eqnarray*}
u_1(x,t)&=&\Psi^*_1 u_0(x,t)+h \mathcal{R}_1\left(\vec{u}_0(x,t)\right)\\
&=& h J^\alpha\left(\mathcal{D}_t^\alpha u_0 -\mu u_0u_{0x}+18u_0u_{0xx}u_{0xx}+36u_{0x}u_{0x}u_{0xx}\right)\\
&&+h J^\alpha\left(24u_0u_{0x}u_{0xxx}+3u_0u_{0}u_{0xxxx}-u_{0xxxx}-u_{0xxxxxx}\right)\\
&=& -\frac{\lambda{h}e^{\lambda x}(\lambda^5+\lambda^3(1-81e^{2\lambda x})+\mu e^{\lambda x}) }{\Gamma(\alpha+1)}t^\alpha
\end{eqnarray*}

\begin{eqnarray*}
u_2(x,t)&=&\Psi^*_2 u_1(x,t)+h J^\alpha\left(\mathcal{R}_2\left(\vec{u}_1(x,t)\right)\right)\\
&=&hu_1+ h J^\alpha\left(\mathcal{D}_t^\alpha u_1-\mu u_{0}u_{1x}-\mu u_{1}u_{0x}+18u_0u_{0xx}u_{1xx}+18u_0u_{1xx}u_{0xx}+18u_1u_{0xx}u_{0xx}\right)\\
&&+ \ h J^\alpha\left(36u_{0x}u_{0x}u_{1xx}+36u_{0x}u_{1x}u_{0xx}+36u_{1x}u_{0x}u_{0xx}+24u_0u_{0x}u_{1xxx}+24u_0u_{1x}u_{0xxx}\right)\\
&&+ \ h J^\alpha\left(24u_{1}u_{0x}u_{0xxx}+3u_{0}u_{0}u_{1xxxx}+3u_{0}u_{1}u_{0xxxx}+3u_{1}u_{0}u_{0xxxx}-u_{1xxxx}-u_{1xxxxxx}\right)\\
&=& (n+h)\,u_1+\frac{\lambda^2h^2e^{\lambda x}\left(151875\lambda^6 e^{4\lambda x}-1092\mu \lambda^3e^{3\lambda x}-3(324\lambda^6(61\lambda^2+7)-\mu^2)e^{2\lambda x}\right)}{\Gamma(2\alpha+1)}t^{2\alpha}\\
&&+\frac{\lambda^2h^2e^{\lambda x}\left((\lambda^2+1)^2\lambda^6+6\mu(11\lambda^2+3) \lambda^3e^{\lambda x}\right)}{\Gamma(2\alpha+1)}t^{2\alpha},
\end{eqnarray*}

\begin{eqnarray*}
u_3(x,t)&=&\Psi^*_3 u_2(x,t)+h J^\alpha\left(\mathcal{R}_3\left(\vec{u}_2(x,t)\right)\right)\\
&=&hu_2+ h J^\alpha\left(\mathcal{D}_t^\alpha u_2-\mu u_{0}u_{2x}-\mu u_{1}u_{1x}-\mu u_{2}u_{0x}+18u_0u_{0xx}u_{2xx}+18u_0u_{1xx}u_{1xx}+18u_1u_{0xx}u_{1xx}\right)\\
&&+ \ h J^\alpha\left(18u_0u_{2xx}u_{0xx}+18u_1u_{1xx}u_{0xx}+18u_2u_{0xx}u_{0xx}+36u_{0x}u_{0x}u_{2xx}+36u_{0x}u_{1x}u_{1xx}+36u_{1x}u_{0x}u_{1xx}\right)\\
&&+ \ h J^\alpha\left(36u_{0x}u_{2x}u_{0xx}+36u_{1x}u_{0x}u_{1xx}+36u_{2x}u_{0x}u_{0xx}+24u_0u_{0x}u_{2xxx}+24u_0u_{1x}u_{1xxx}+24u_1u_{0x}u_{1xxx}\right)\\
&&+ \ h J^\alpha\left(24u_0u_{2x}u_{0xxx}+24u_0u_{1x}u_{1xxx}+24u_2u_{0x}u_{0xxx}+3u_{0}u_{0}u_{2xxxx}+3u_{0}u_{1}u_{1xxxx}+3u_{1}u_{0}u_{1xxxx}\right)\\
&&+ \ h J^\alpha\left(3u_{0}u_{2}u_{0xxxx}+3u_{1}u_{1}u_{0xxxx}+3u_{2}u_{0}u_{0xxxx}-u_{2xxxx}-u_{2xxxxxx}\right)\\
&=& (n+h)\,u_2+\frac{\lambda^2h^2(n+h)e^{\lambda x}\left((\lambda^2+1)^2\lambda^6+6\mu(11\lambda^2+3)\lambda^3e^{\lambda x}\right)}{\Gamma(2\alpha+1)}t^{2\alpha}\\
&&+\frac{\lambda^2h^2(n+h)e^{\lambda x}\left(151875\lambda^6 e^{4\lambda x}-1092\mu \lambda^3e^{3\lambda x}-3(324\lambda^6(61\lambda^2+7)-\mu^2)e^{2\lambda x}\right)}{\Gamma(2\alpha+1)}t^{2\alpha}\\
&&+\frac{\lambda^{12}h^3e^{\lambda x}\left(-151875(16357\lambda^2+709)e^{4\lambda x}-(\lambda^2+1)^3+1093955625e^{6\lambda x}\right)}{\Gamma(3\alpha+1)}t^{3\alpha}\\
&&+\frac{\lambda^{12}h^3e^{\lambda x}\left(243(177877\lambda^4+40178\lambda^2+2269)e^{6\lambda x}\right)}{\Gamma(3\alpha+1)}t^{3\alpha}\\
&&+\frac{3\mu\lambda^6h^3e^{3\lambda x} \left(1586896\lambda^5e^{\lambda x}-1718982\lambda^3e^{3\lambda x}-795\mu \lambda^2+3695\mu e^{2\lambda x} \right)}{\Gamma(3\alpha+1)}t^{3\alpha}\\
&&+\frac{\mu \lambda^3h^3e^{2\lambda x}\left(12(26716\lambda^6-\mu^2)e^{2\lambda x}-297\mu \lambda^3e^{\lambda x}-2(2113\lambda^4+1106\lambda^2+145)\lambda^6\right)}{\Gamma(3\alpha+1)}t^{3\alpha}\\
&&+\frac{h^3e^{3\lambda x}\Gamma(2\alpha+1)\left(243\lambda^{16}+\lambda^{14}(486-303750e^{2\lambda x})+\lambda^{12}(47258883e^{4\lambda x}-303750e^{2\lambda x}+243)-3\mu^2 \lambda^6\right)}{\Gamma(\alpha+1)^2\Gamma(3\alpha+1)}t^{3\alpha}\\
&&-\frac{\mu h^3e^{2\lambda x}\Gamma(2\alpha+1)\left(\lambda^{10}-1860\lambda^8e^{2\lambda x}2\lambda^{8}-1860\lambda^{6}e^{2\lambda x}+649539\lambda^{6}e^{4\lambda x}+\lambda^6\right)}{\Gamma(\alpha+1)^2\Gamma(3\alpha+1)}t^{3\alpha}\\
&&-\frac{\mu h^3e^{2\lambda x}\Gamma(2\alpha+1)\left(3\mu \lambda^5e^{\lambda x}-2280\mu \lambda^3e^{3\lambda x}+2\mu^2 e^{2\lambda x}\right)}{\Gamma(\alpha+1)^2\Gamma(3\alpha+1)}t^{3\alpha}.
\end{eqnarray*}
Using the same procedure, expression for $u_m(x,t),$ $m=4,5,6,...$ can be obtained. The expression of the series solutions given by q-HAM can be written in the form

\begin{eqnarray}\notag 
u(x,t;n;h)&\cong& u_0(x,t)+\sum_{i=1}^3 u_i(x,t;n;h)\left(\frac{1}{n}\right)^i \\ \notag
&=&e^{\lambda x}-\frac{\lambda{h}(3n^2+3nh+h^2)e^{\lambda x}(\lambda^5+\lambda^3(1-81e^{2\lambda x})+\mu e^{\lambda x}) }{n^3\Gamma(\alpha+1)}t^\alpha\\ \notag
&&+\frac{\lambda^2h^2(3n+2h)e^{\lambda x}\left(151875\lambda^6 e^{4\lambda x}-1092\mu \lambda^3e^{3\lambda x}-3(324\lambda^6(61\lambda^2+7)-\mu^2)e^{2\lambda x}\right)}{n^3\Gamma(2\alpha+1)}t^{2\alpha}\\ \notag
&&+\frac{\lambda^2h^2(3n+2h)e^{\lambda x}\left((\lambda^2+1)^2\lambda^6+6\mu(11\lambda^2+3) \lambda^3e^{\lambda x}\right)}{n^3\Gamma(2\alpha+1)}t^{2\alpha}
\\ \notag
&&+\frac{\lambda^{12}h^3e^{\lambda x}\left(-151875(16357\lambda^2+709)e^{4\lambda x}-(\lambda^2+1)^3+1093955625e^{6\lambda x}\right)}{n^3\Gamma(3\alpha+1)}t^{3\alpha}\\ \notag
&&+\frac{\lambda^{12}h^3e^{\lambda x}\left(243(177877\lambda^4+40178\lambda^2+2269)e^{6\lambda x}\right)}{n^3\Gamma(3\alpha+1)}t^{3\alpha}\\ \notag
&&+\frac{3\mu\lambda^6h^3e^{3\lambda x} \left(1586896\lambda^5e^{\lambda x}-1718982\lambda^3e^{3\lambda x}-795\mu \lambda^2+3695\mu e^{2\lambda x} \right)}{n^3\Gamma(3\alpha+1)}t^{3\alpha}\\ \notag
&&+\frac{\mu \lambda^3h^3e^{2\lambda x}\left(12(26716\lambda^6-\mu^2)e^{2\lambda x}-297\mu \lambda^3e^{\lambda x}-2(2113\lambda^4+1106\lambda^2+145)\lambda^6\right)}{n^3\Gamma(3\alpha+1)}t^{3\alpha}\\ \notag
&&+\frac{h^3e^{3\lambda x}\Gamma(2\alpha+1)\left(243\lambda^{16}+\lambda^{14}(486-303750e^{2\lambda x})+\lambda^{12}(47258883e^{4\lambda x}-303750e^{2\lambda x}+243)-3\mu^2 \lambda^6\right)}{n^3\Gamma(\alpha+1)^2\Gamma(3\alpha+1)}t^{3\alpha}\\ \notag
&&-\frac{\mu h^3e^{2\lambda x}\Gamma(2\alpha+1)\left(\lambda^{10}-1860\lambda^8e^{2\lambda x}2\lambda^{8}-1860\lambda^{6}e^{2\lambda x}+649539\lambda^{6}e^{4\lambda x}+\lambda^6\right)}{n^3\Gamma(\alpha+1)^2\Gamma(3\alpha+1)}t^{3\alpha}\\ 
&&-\frac{\mu h^3e^{2\lambda x}\Gamma(2\alpha+1)\left(3\mu \lambda^5e^{\lambda x}-2280\mu \lambda^3e^{3\lambda x}+2\mu^2 e^{2\lambda x}\right)}{n^3\Gamma(\alpha+1)^2\Gamma(3\alpha+1)}t^{3\alpha}. \label{eq:solut5cahnqham4}
\end{eqnarray}
Thus, \autoref{eq:solut5cahnqham4} give an approximate solution to the IVP \ref{eq:cahn2}-\eqref{eq:cahn2IC} in terms of convergence parameters $h$ and $n$.
%%%%%%%%%%%%%%%%%%%%%%%%%%%%%%%%%%%%%%%%%%%%%%%%%%%%%%%%%%%%%%%%%
        %%%%%%%%%% Numerical Results %%%%%%%%%%%%%%
%%%%%%%%%%%%%%%%%%%%%%%%%%%%%%%%%%%%%%%%%%%%%%%%%%%%%%%%%%%%%%%%%
\subsection{Numerical Results for TFCH Equation of Sixth-Order}
\begin{figure}[H]
    \centering
    \begin{subfigure}[b]{0.40\textwidth}
        \includegraphics[width=\textwidth]{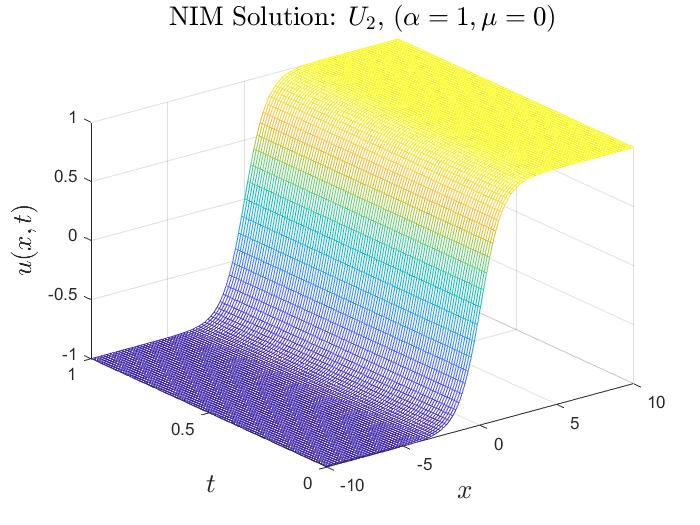}
        %\caption{}
        %\label{Fig1u}
    \end{subfigure}
    ~ %add desired spacing between images, e. g. ~, \quad, \qquad, \hfill etc. 
    \begin{subfigure}[b]{0.40\textwidth}
        \includegraphics[width=\textwidth]{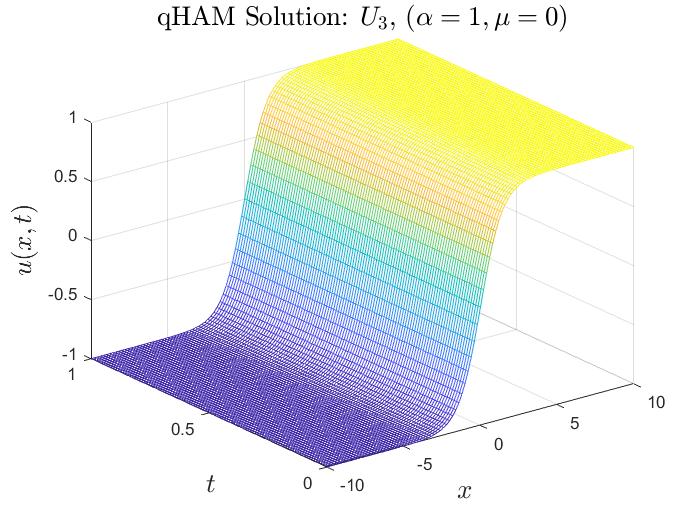}
        %\caption{}
        %\label{fi2}
\end{subfigure}
\caption{Case I: Comparison between NIM and q-HAM Solution for $n=1,$ and $h=-1.$}\label{fig:Fig8}
\end{figure}

\begin{figure}[H]
    \centering
    \begin{subfigure}[b]{0.40\textwidth}
        \includegraphics[width=\textwidth]{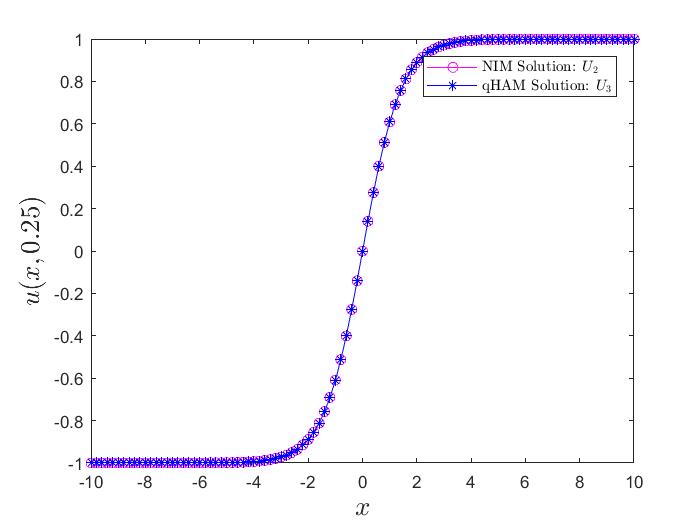}
        %\caption{}
        %\label{Fig1u}
    \end{subfigure}
\caption{Case I: Comparison between NIM and q-HAM Solution in 2D for $n=\alpha=1,\ \mu=0,$ and $h=-1.$}\label{fig:Fig9}
\end{figure}

This section present the numerical simulation of each case to demonstrate the effectiveness of the two iterative methods used for solving TFCH \autoref{eq:1N} subject to different initial conditions. \Cref{fig:Fig8,fig:Fig9} show the NIM $U_2$- and q-HAM $U_3$- solutions for $\alpha=1$ and $\mu=0$ which are graphically almost indistinguishable.
%%%%%%%%%%%%%%%%%%%%%%%%%%%%%%%%%%%%%%%%%%%%%%%%%%%
%%%%%%%%%%%%%%%% 4 Case III Graph %%%%%%%%%%%%%%%%
%%%%%%%%%%%%%%%%%%%%%%%%%%%%%%%%%%%%%%%%%%%%%%%%%%%%
%%%%%%%%%%%%%%%%%%%%%%%%%%%%%%%%%%%%%%%%%%%%%%%%%%%%
%%%%%%%%%%%%%%%%%%%%%%%%%%%%%%%%%%%%%%%%%%%%%%%%%%%%%
%%%%%%%%%%%%%%%%%%%%%%%%%%%%%%%%%%%%%%%%%%%%%%%%%%%%%
\begin{figure}[H]
    \centering
    \begin{subfigure}[b]{0.40\textwidth}
\includegraphics[width=\textwidth]{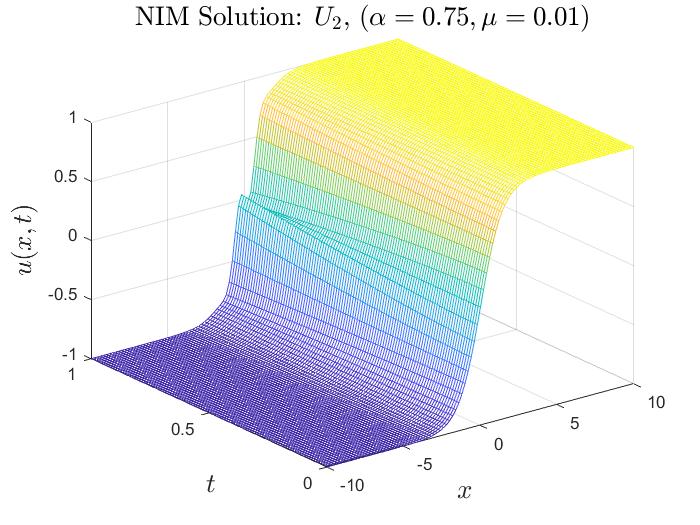}
        %\caption{}
        %\label{Fig1u}
    \end{subfigure}
    ~ %add desired spacing between images, e. g. ~, \quad, \qquad, \hfill etc. 
    \begin{subfigure}[b]{0.40\textwidth}
\includegraphics[width=\textwidth]{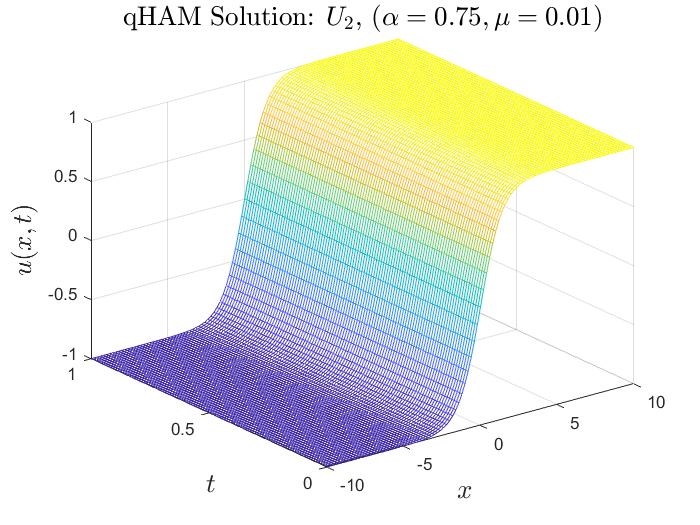}
        %\caption{}
        %\label{fi2}
    \end{subfigure}
\caption{Case I: Comparison between NIM and q-HAM Solution for $n=20$ and $h=-0.4.$}\label{fig:Fig10}
\end{figure} 
%%%%%%%%%%%%%%%%%%%%%%%%%%%%%%%%%%%%%%%%%%%%%%%%%%%%%
\begin{figure}[H]
    \centering
    \begin{subfigure}[b]{0.40\textwidth}
\includegraphics[width=\textwidth]{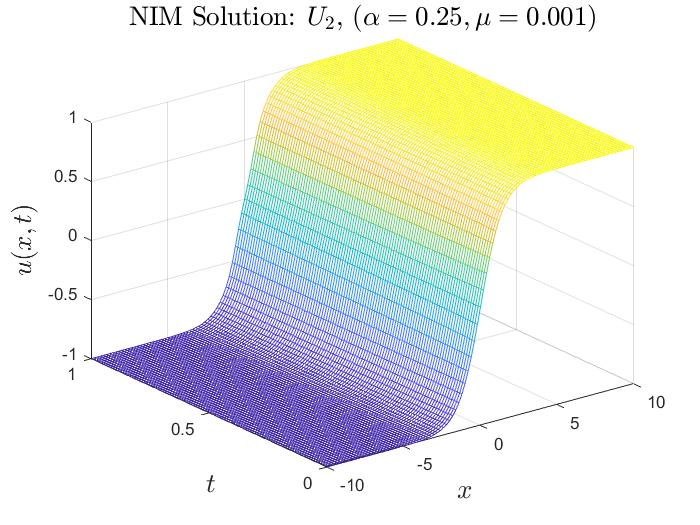}
        %\caption{}
        %\label{Fig1u}
    \end{subfigure}
    ~ %add desired spacing between images, e. g. ~, \quad, \qquad, \hfill etc. 
    \begin{subfigure}[b]{0.40\textwidth}
\includegraphics[width=\textwidth]{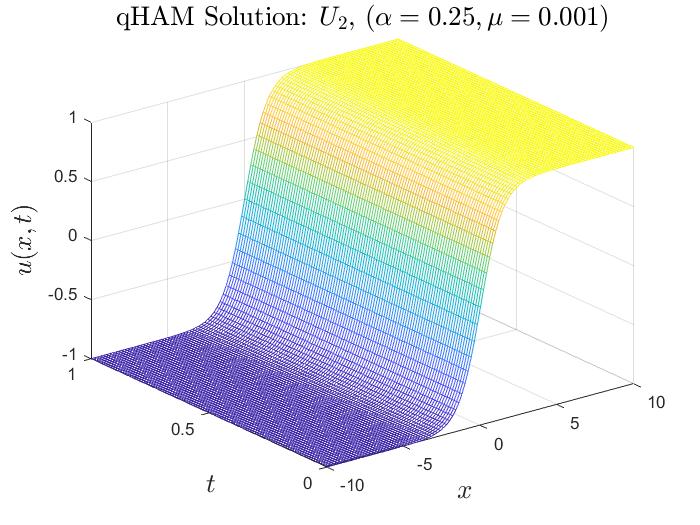}
        %\caption{}
        %\label{fi2}
    \end{subfigure}
\caption{Case I: Comparison between NIM and q-HAM Solution for $n=20$ and $h=-0.4.$}\label{fig:Fig11}
\end{figure}  
%%%%%%%%%%%%%%%%%%%%%%%%%%%%%%%%%%%%%%%%%%%%%%%%%%%
     %%%%%%%%% Case II Graph %%%%%%%%%%%%%%%%
%%%%%%%%%%%%%%%%%%%%%%%%%%%%%%%%%%%%%%%%%%%%%%%%%%%
\begin{figure}[H]
    \centering
    \begin{subfigure}[b]{0.40\textwidth}
\includegraphics[width=\textwidth]{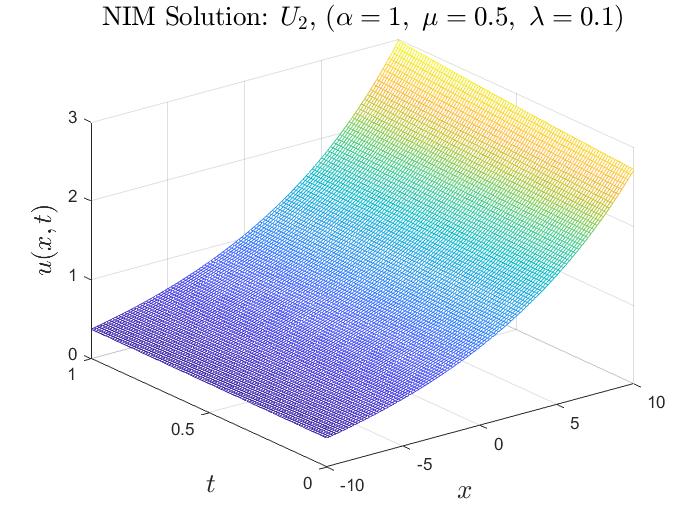}
        %\caption{}
        %\label{Fig1u}
    \end{subfigure}
    ~ %add desired spacing between images, e. g. ~, \quad, \qquad, \hfill etc. 
    \begin{subfigure}[b]{0.40\textwidth}
    \includegraphics[width=\textwidth]{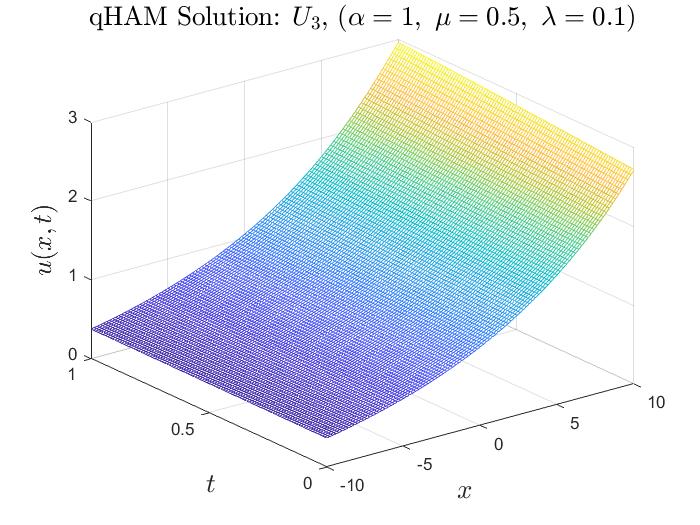}
        %\caption{}
        %\label{fi2}
    \end{subfigure}
    \caption{Case II: Comparison between NIM and q-HAM Solution for $n=1$ and $h=-0.8.$}\label{fig:Fig12}
\end{figure} 
%%%%%%%%%%%%%%%%%%%%%%%%%%%%%%%%%%%%%%%%%%%%%%%%%%%
\begin{figure}[H]
    \centering
    \begin{subfigure}[b]{0.40\textwidth}
\includegraphics[width=\textwidth]{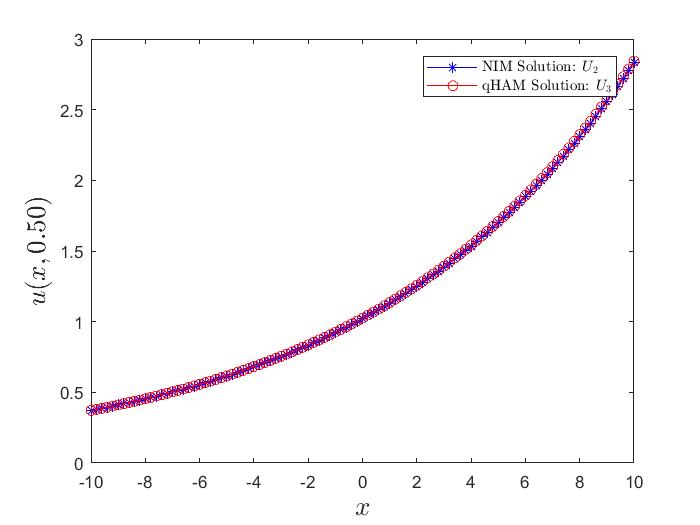}
        %\caption{}
        %\label{Fig1u}
    \end{subfigure}
\caption{Case II: Comparison between NIM and q-HAM Solution in 2D for $n=\alpha=1,\ \mu=0.5,$\ $h=-0.8,$ and $\lambda=0.1.$}\label{fig:Fig13}
\end{figure} 
%%%%%%%%%%%%%%%%%%%%%%%%%%%%%%%%%%%%%%%%%%%%%%%%%%%%%%%%%%%%%%%%%%%%
\begin{figure}[H]
    \centering
    \begin{subfigure}[b]{0.40\textwidth}
    \includegraphics[width=\textwidth]{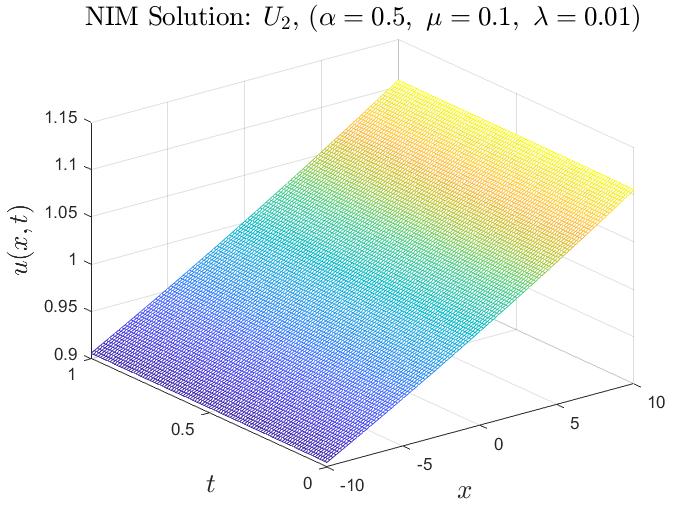}
        %\caption{}
        %\label{Fig1u}
    \end{subfigure}
    ~ %add desired spacing between images, e. g. ~, \quad, \qquad, \hfill etc. 
    \begin{subfigure}[b]{0.40\textwidth}
    \includegraphics[width=\textwidth]{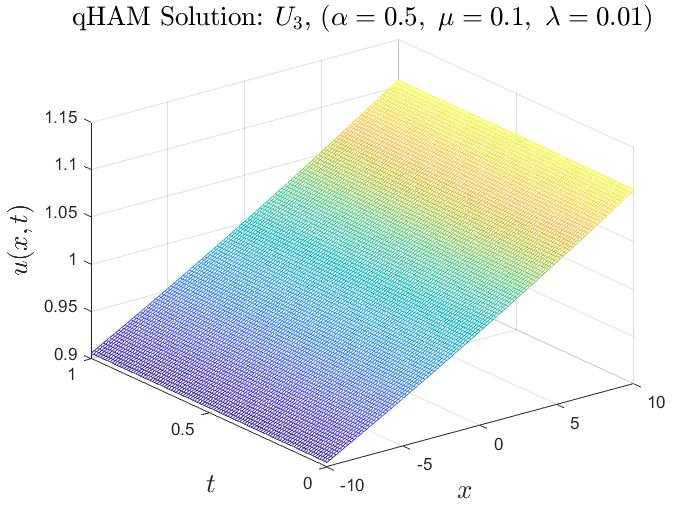}
        %\caption{}
        %\label{fi2}
    \end{subfigure}
\caption{Case II: Comparison between NIM and q-HAM Solution for $n=1$ and $h=-1$.}\label{fig:Fig14}
\end{figure} 
%%%%%%%%%%%%%%%%%%%%%%%%%%%%%%%%%%%%%%%%%%%%%%%%%%%%%%%%%%%%%%%%%%%%
\begin{figure}[H]
    \centering
    \begin{subfigure}[b]{0.40\textwidth}
    \includegraphics[width=\textwidth]{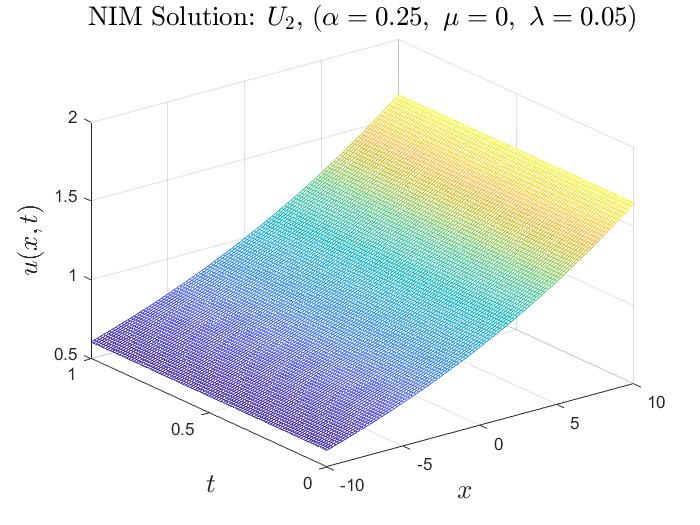}
        %\caption{}
        %\label{Fig1u}
    \end{subfigure}
    ~ %add desired spacing between images, e. g. ~, \quad, \qquad, \hfill etc. 
    \begin{subfigure}[b]{0.40\textwidth}
    \includegraphics[width=\textwidth]{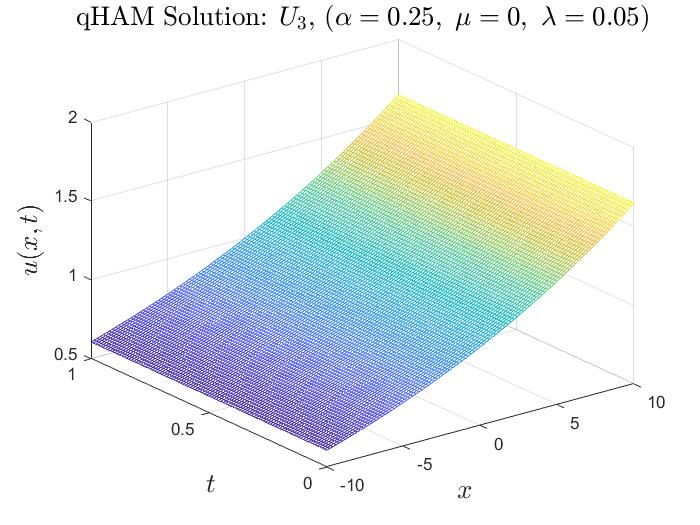}
        %\caption{}
        %\label{fi2}
    \end{subfigure}
\caption{Case II: Comparison between NIM and q-HAM Solution for $n=1$ and $h=-1$.}\label{fig:Fig15}
\end{figure} 
\begin{remark}
In \Crefrange{fig:Fig8}{fig:Fig15}, we observed that the fraction factor $\Big(\frac{1}{n}\Big)^i$ and the parameter $h$ highly increase the convergence of the chances of the q-HAM.
\end{remark}
%%%%%%%%%%%%%%%%%%%%%%%%%%%%%%%%%%%%%%%%%%%%%%%%%%%%%%%%%%%%%%%%%
         %%%%%%%%%%%%%% Conclusion %%%%%%%%%%%%%%%%
%%%%%%%%%%%%%%%%%%%%%%%%%%%%%%%%%%%%%%%%%%%%%%%%%%%%%%%%%%%%%%%%%
%\newpage
\section{Concluding remarks}
In conclusion, we have studied iterative methods of constructing approximate solutions to the time-fractional nonlinear Cahn-Hilliard equations of fourth and sixth-order with different initial conditions. We used NIM and q-HAM to obtain approximate series solutions. As shown in our examples, the two iterative methods does not required any transformations or perturbation. Therefore, these methods are considerably efficient, powerful, and easy to implement when compared to other numerical methods in constructing approximated solutions to the linear and nonlinear differential equations. Our aim in this paper is not to conclude that one method is better than the other, but rather, conclude that both methods provides a good approximate solutions, even is some cases we can obtain the exact solutions.
%%%%%%%%%%%%%%%%%%%%%%%%%%%%%%%%%%%%%%%%%%%%%%%%%%%%%%%%%%%%%%%%%
           %%%%%%%%%%% bibliography %%%%%%%%%%%%%%%%%
%%%%%%%%%%%%%%%%%%%%%%%%%%%%%%%%%%%%%%%%%%%%%%%%%%%%%%%%%%%%%%%%%
%\newpage


\begin{thebibliography}{99}
\bibitem{Leib} Leibniz GW. Letter from Hanover, Germany to G.F.A. LíHospital, September 30, 1695, Leibniz Mathematische Schriften, Olms-Verlag, Hildesheim, Germany. 1962;301-302. First published in 1849.

\bibitem{Leib1} Leibniz GW. Letter from Hanover, Germany to Johann Bernoulli, December 28, 1695, Leibniz Mathematische Schriften, Olms-Verlag, Hildesheim, Germany. 1962;226. First published in 1849.

\bibitem{fc} Podlubny I. Fractional differential equations: an introduction to fractional derivatives, fractional differential equations, to methods of their solution and some of their applications. New York: Academic Press; 1999.

\bibitem{ChenChen} Chen D, Chen Y, Xue D. Three fractional-order TV-models for image de-noising. J Comput Inf Sys. 2013;\textbf{9}(12):4773-80.

\bibitem{Hilf1} Hilfer R, Anton L. Fractional master equations and fractal time random walks. Physical Review E. 1995;51:R848-R851.

\bibitem{Main} Mainardi F. Fractional Calculus and Waves in Linear Viscoelasticity. Imperial College Press; 2010.

\bibitem{Monj} Monje CA, Chen Y, Vinagre BM, Xue D, Feliu V. Fractional-order Systems and Controls, Advances in Industrial Control. Springer; 2010.

\bibitem{PU} Pu YF. Fractional differential analysis for texture of digital image. J Alg Comput Technol. 2007;\textbf{1}(03):357-80.

\bibitem{Bal} Sun HG, Zhang Y, Baleanu D, Chen W, Chen YQ. A new collection of real world applications of fractional calculus in science and engineering.  Commun. in Nonlinear Sci. Numer. Simulat. 2018;\textbf{64}:213-231.

\bibitem{VEVV} Tarasov VE, Tarasova VV. Time-dependent fractional dynamics with memory in quantum and economic physics. Ann. Phys. 2017:383,579-599.

\bibitem{Ucs} Ullah A, Chen W, Sun HG,  Khan MA. An efficient variational method for restoring images with combined additive and multiplicative noise. Int J Appl Comput Math. 2017;\textbf{3}(3):1999-2019.

\bibitem{ZW} Zhang J, Wei Z, Xiao L.  Adaptive fractional multiscale method for image de-noising. J Math Imaging Vis. 2012;43:39-49.

\bibitem{ZPZ} Zhang Y, Pu YF, Hu JR, Zhou JL. A class of fractional-order variational image in-painting models. Appl Math Inf Sci. 2012;\textbf{06}(02):299-306.

\bibitem{Y} Zhang Y, Benson DA, Meerschaert MM, LaBolle EM, Scheffler  HP. Random walk approximation of fractionalorder multiscaling anomalous diffusion. Physical Review E. 2006;74:026706-026715.

\bibitem{CahnE} Cahn JW, Hilliard JE. “Free energy of a non-uniform system I: Interfacial free energy.” The Journal of Chemical Physics. 1958:28(2):258-267.

\bibitem{Cahn} Cahn JW. On spinodal decomposition. Acta. Metall.. 1961;9:795-801.

\bibitem{EllZ0} Elliott CM. The Cahn-Hilliard model for the kinetics of phase separation, in Mathematical Models for Phase Change Problems, (ed. J. F. Rodrigues). International Series of Numerical Mathematics, 88, Birkh¨auser, Basel; 1989.

\bibitem{NOV} Novick-Cohen A. “The Cahn-Hilliard equation,” in Handbook of Differential Equations, Evolutionary Partial Differential Equations, 4, (eds. C. M. Dafermos and M. Pokorny). Elsevier, Amsterdam. 2008;201-228.

\bibitem{MM} Mohamed M, Mekheimer K. Analytical approximate solution for nonlinear space–time fractional Cahn–Hilliard equation. Int. Electron. J. Pure Appl. Math. 2014;7:145-149.

\bibitem{EllZ1} Copetti MIM, Elliott CM. Kinetics of phase decomposition processes: numerical solutions to the Cahn—Hilliard equation. Materials Science and Technology. 1990;6:273-283.

\bibitem{DahBen} Dhmani Z, Benbachir M. “Solutions of the Cahn-Hilliard equation with time- and space-fractional derivatives.” International Journal of Nonlinear Science. 2009;8(1):19-26.

\bibitem{DeM} De Mello EVL, Otto T, Da Silveira F. Numerical study of the Cahn-Hilliard equation in one, two and three dimensions. Physica A. 2005;347: 429-443.

\bibitem{Jun} Junseok K. A numerical method for the Cahn-Hilliard equation with a variable mobility. Commun. Nonlinear Sci. Numer. Simul. 2007,12: 1560-1572.

\bibitem{RyHoff} Rybka  P, Hoffmann KH, Convergence of solutions to Cahn-Hilliard equation, Comm. Partial Differential Equations. 1999;24:1055-1077.

\bibitem{UK} Ugurlu Y, Kaya D. Solutions of the Cahn-Hilliard equation. Comput. Math. Appl. 2008;56:3038-3045.

\bibitem{adm2} Ray SS, Bera RK. Analytical solution of a fractional diffusion equation by Adomian decomposition method. Applied Mathematics and Computation. 2006;174 (1):329-336.

\bibitem{Das} Das S. Analytical solution of a fractional diffusion equation by variational iteration method, Comput. Math. Appl. 2009;57:483-487.

\bibitem{vim} Neamaty A, Agheli B, Darzi R. Variational iteration method and He's polynomials for time-fractional partial differential equations. Progress in Fractional Differentiation and Applications. 2015;1(1):47-55.

\bibitem{Atangana1} Atangana A, Alabaraoye E. Solving system of fractional partial differential equations arisen in the model of HIV infection of CD4+ cells and attractor one-dimensional Keller- Segel equation. Advances in Difference Equations. 2013;2013:94.

\bibitem{Atangana2} Atangana A,  Kilicman A. Analytical solutions of the space time fractional derivative of advection dispersion equation. Mathematical Problems in Engineering, vol. 2013, Article ID 853127, 9 pages, 2013.  https://doi.org/10.1155/2013/853127.

\bibitem{khals} Khalid M, Sultana M, Zaidi  F, Uroosa A. “Solving Linear and Nonlinear Klein-Gordon Equations by New Perturbation Iteration Transform Method.” TWMS J. App. Eng. Math. 2016;\textbf{6}(1):115-125.

\bibitem{hpm} Abdulaziz O, Hashim I, Momani S. Application of homotopy-perturbation method to fractional IVPs. Journal of Computational and Applied Mathematics. 2008;\textbf{216}(2);574-584.

\bibitem{Alquran2} Alquran, M. Analytical solution of time-fractional two-component evolutionary system of order 2 by residual power series method. J. Appl. Anal. Comput. 2015; 589-599.

\bibitem{Xu} Xu F, Gao Y, Yang X, Zhang H. Construction of fractional power series solutions to fractional Boussinesq equations using residual power series method, Math. Probl. Eng. 2016;(2016): 15. Article ID 5492535

\bibitem{Bha} Bhalekar S. Daftardar-Gejji V. New iterative method: application to partial differential equations. Applied Mathematics and Computation. 2008;\textbf{203}(2):778-783.

\bibitem{Daf1} Daftardar-Gejji V, Bhalekar S. Solving fractional boundary value problems with Dirichlet boundary conditions using a new
iterative method. Computers \& Mathematics with Applications. 2010;\textbf{59}(5):1801-1809. 

\bibitem{Daf} Daftardar-Gejji V, Jafari H. An iterative method for solving nonlinear functional equations. Journal of Mathematical Analysis and Applications. 2006;\textbf{316}(2):753-763.

\bibitem{Wen} Wenjin Li, Yanni Pang. An Iterative Method for Time-Fractional Swift-Hohenberg Equation.  Advances in Mathematical Physics. 2018;(1-2):1-13.

\bibitem{Taw} El-Tawil MA, Huseen SN. The Q-homotopy analysis method (Q-HAM). Int. J. Appl. Math. Mech. 2012;\textbf{8}(15):51-75.

\bibitem{IyizOla} Iyiola OS, Exact and Approximate Solutions of Fractional Diffusion Equations with Fractional Reaction Terms. Progr. Fract. Differ. Appl. 2016;\textbf{2}(1):21-30.
 
\bibitem{Iyiz1} Iyiola OS, Zaman FD. A fractional diffusion equation model for cancer tumor. AIP Adv. 2014;4:107-121.

\bibitem{Iyiz3} Iyiola OS, Soh ME, Enyi CD. Generalised homotopy analysis method (q-HAM) for solving foam drainage equation of time fractional type. Math. Eng. Sci. Aerospace. 2013;\textbf{4}(4):105.

\bibitem{Miller} Miller KS, Ross B. An introduction to the fractional calculus and fractional differential equations. John Wiley and Sons, New York; 2003.

\bibitem{kil} Kilbas AA, Srivastava HM, Trujillo JJ. Theory and Applications of Fractional Differential Equations, North-Holland Mathematical Studies, Vol. 204, Elsevier (North-Holland) Science Publishers, Amsterdam, London and New York; 2006.

\bibitem{Liao} Liao SJ. An approximate solution technique not depending on small parameters: a special example. Int. J. Nonlin. Mech. 1995;\textbf{30}(3):371-380.
\end{thebibliography}
\end{document}